\documentclass[12pt]{amsart}  
\usepackage[T1]{fontenc}
\usepackage{amsxtra}

\usepackage{color}
\usepackage{amscd,amsfonts,amsmath,amssymb,
amsthm,amstext,latexsym,amstext}
\usepackage{enumerate,pifont,graphicx}
\usepackage[english]{babel}
\usepackage[all,cmtip]{xy} % for commutative diagram
\newcommand{\R}{{\mathbb R}}
\newcommand{\C}{{\mathbb C}}
\newcommand{\N}{{\mathbb N}}
\renewcommand{\S}{{\mathbb S}}
\renewcommand{\H}{{\mathbb H}}
\renewcommand{\L}{{\mathbb L}}
\def\bfa{{\bf a}}
\def\bfb{{\bf b}}
\def\bfq{{\bf q}}
\def\bft{{\bf t}}
\def\bfp{{\bf p}}
\def\[{{[\![}}
\def\]{{]\!]}}

\def\boF{{\mathcal F}}

\def\boM{{\mathcal M}}

\def\CC{\overline{\C}}
\def\cqfd{\hfill$\Box$}
\def\Res{{\,\rm Res}}

\def\Re{{\rm Re}}
\def\Im{{\rm Im}}

\def\i{{\rm i}}

\def\wtx{{\widetilde{x}}}

\def\wtA{{\widetilde{A}}}

\def\wtF{{\widetilde{F}}}
\def\wtG{{\widetilde{G}}}

\def\wtY{{\widetilde{Y}}}

\def\myomega{{\Omega}}
\def\Lambda{{M}}
\def\whLambda{{\widehat{M}}}
\def\whlambda{{\widehat{\lambda}}}
\def\whPi{{\widehat{\Pi}}}
\def\whA{{\widehat{A}}}
\def\whF{{\widehat{F}}}
\def\whf{{\widehat{f}}}
\def\whx{{\widehat{x}}}
\def\whG{{\widehat{G}}}

\def\whS{{\widehat{S}}}
\def\half{\mbox{$\frac{1}{2}$}}

\newtheorem{theorem}{Theorem}
\newtheorem{question}{Question}
\newtheorem{lemma}{Lemma}
\newtheorem{proposition}{Proposition}
\newtheorem{remark}{Remark}
\newtheorem{corollary}{Corollary}

\newtheorem{definition}{Definition}

\title{Opening nodes on horosphere packings}
\author{Martin Traizet}
\begin{document}
\maketitle
\section{Introduction}
In this paper, a horosphere packing will be a finite, connected set of distinct
horospheres in hyperbolic space $\H^3$, such that any two horospheres are
either disjoint or tangent. We will denote $n$ the number of horospheres
and $m$ the number of tangency points.

\medskip

Horospheres have constant mean curvature equal to $1$ (CMC-1 for short).
In \cite{pacard}, Pacard and Pimentel have constructed complete, embedded CMC-1
surfaces in hyperbolic space by desingularization of a horosphere packing.
Heuristically, the horospheres are slightly shrinked and a suitable small catenoid
is glued at each tangency point (see Figure 1).

\medskip

On the other hand, CMC-1 surfaces in hyperbolic space
admit a Weierstrass-type representation in term of meromorphic data, as
discovered by Bryant \cite{bryant}. In particular,
they have a meromorphic Gauss map $G$ and a holomorphic quadratic
differential $Q$.
Our goal in this paper is to revisit the result of Pacard and Pimentel using
Bryant's representation. We prove:
\begin{theorem}
\label{main-theorem}
Given a packing of $n$ horospheres with $m$ tangency points, there exists a smooth family
$(M_s)_{0<s<\varepsilon}$ of complete, embedded CMC-1 surfaces in
hyperbolic space $\H^3$, such that
$M_s$ converges when $s\to 0$ to the given horosphere packing.
The surfaces $M_s$ have genus $m-n+1$ and $n$ catenoid-cousin-type ends.
They have finite total curvature.
\end{theorem}

We prove this theorem using Bryant's representation.
Our point of view on Riemann surfaces is that of opening nodes.
The construction follows the general
strategy developped by the author in \cite{traizet,nosym,triply} to construct minimal
surfaces in euclidean space $\R^3$ by opening nodes, using the Weierstrass representation.
The new point is that instead of having to solve a period problem, which is homology
invariant,
we have to solve a monodromy problem, which is homotopy invariant hence
more difficult.
\medskip

The surfaces we construct will in fact depend on $n$ complex parameters
$c_1,\cdots,c_n$, which are the limit points of the horospheres on the ideal boundary
of $\H^3$ (identified with the Riemann sphere) and $n$ positive real parameters $\xi_1,\cdots,\xi_n$
which represent the speed at which each horosphere is ``deflated'' to accomodate the
catenoidal necks, where we think of $s$ as the time parameter. We can normalize
$\xi_1=1$, so our family depends on $3n$ real parameters, which is the expected dimension
for the space of CMC-1 surfaces with $n$ ends.

\medskip

\begin{figure}
\label{fig1}
\begin{center}
\includegraphics[height=40mm]{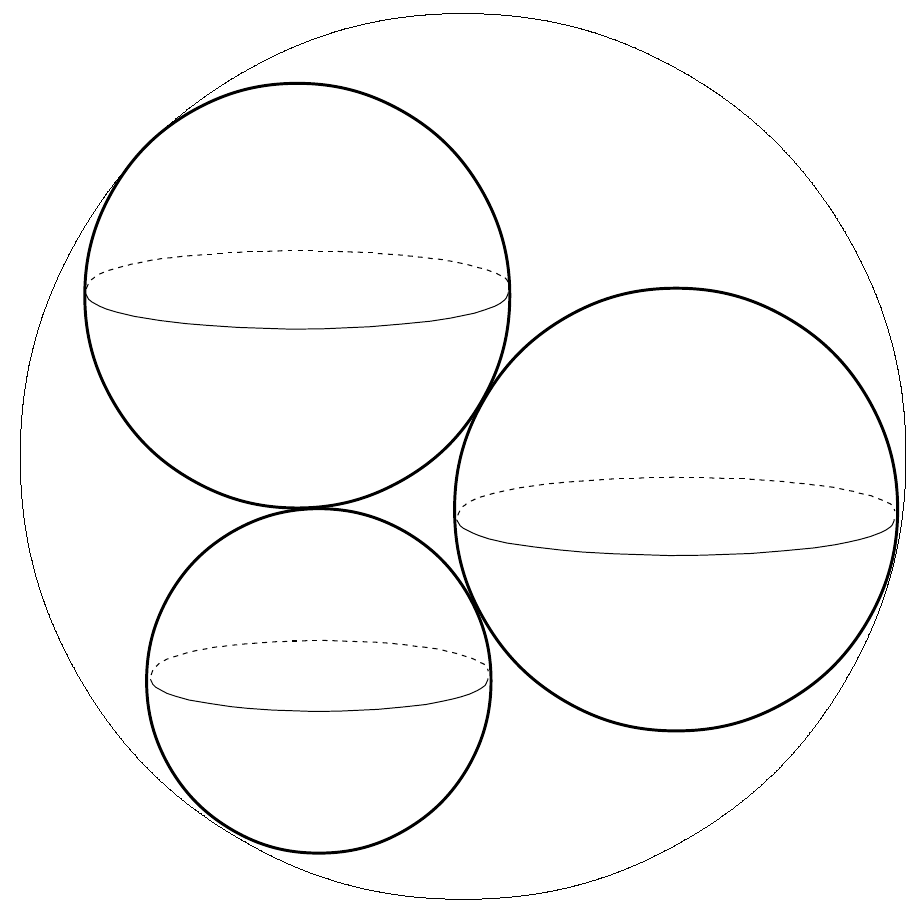}
\hspace{1cm}
\includegraphics[height=40mm]{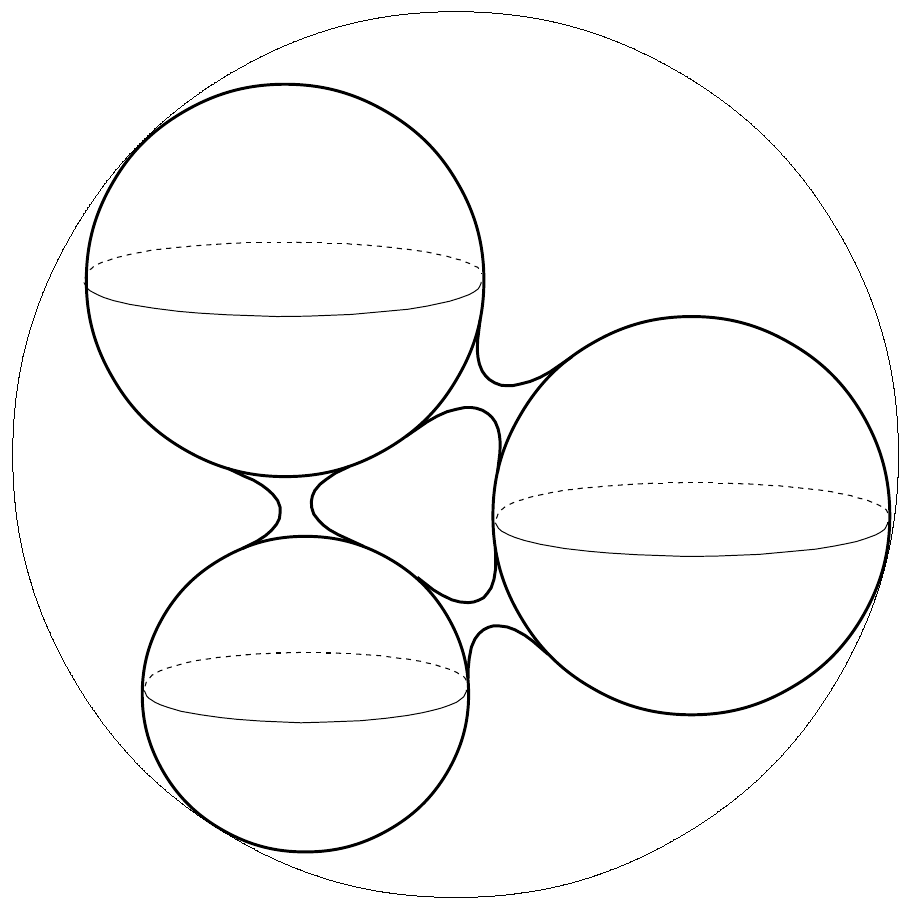}
\end{center}
\caption{Left: a packing of three horospheres with three tangency points
in the ball model of $\H^3$.
Right: a CMC-1 surface of genus one with three ends.
}
\end{figure}

One may ask what are the possible topologies that one can achieve with this
construction. In other words, given the number of ends $n$, what are the possible
genera $g$? Since one can always ignore some tangency points (see Remark
\ref{remark-ignore}), this boils down to the following question:
\begin{question}
\label{question-horosphere}
What is the maximum number of tangency
points that a packing of $n$ horospheres may have?
\end{question}
For $n\geq 3$, Pacard and Pimentel explain how to construct
horosphere packings with up to $3n-6$ tangency points, which gives genus up to $2n-5$.
It turns out that one can do much better!
We investigate Question \ref{question-horosphere} in
Section \ref{section-packing}. We prove that the number of tangency
points is always less than $5n$, and we give examples where it grows like
$4n$. Still, the question remains open.

\medskip

The rest of the paper is organized as follows.
In Section \ref{section-preliminaries},
we explain the principle of Bryant representation of CMC-1 surfaces in
hyperbolic space, and recall standard material about monodromy and
opening nodes.
In Section \ref{section-data}, we construct a family of candidates for the holomorphic
data of the CMC-1 immersions we want to construct.
The monodromy problem is solved in Section \ref{section-monodromy} using the
implicit function theorem.
In Section \ref{section-embedded}, we study the geometry of the surfaces we
have constructed and prove that they are embedded.
Finally, an appendix contains several results of independent interest that we use
in this construction.

\subsection{Related works}
Several authors have used using Bryant representation to study
CMC-1 surfaces in hyperbolic space:
Beno\^it Daniel \cite{daniel1,daniel2,daniel3},
Ricardo Sa Earp and Eric Toubiana \cite{toubiana},
Wayne Rossman, Masaaki Umehara and Kotaro Yamada \cite{rossman},
Masaaki Umehara and Kotaro Yamada \cite{umehara-yamada1,umehara-yamada2,umehara-yamada3}.
Many examples have been constructed which are inspired from known minimal surfaces in $\R^3$.

\medskip

Pascal Collin, Laurent Hauswirth and Harold Rosenberg \cite{collin} have proved the fundamental result
that a properly embedded
CMC-1 surface in $\H^3$ of finite topology must have finite total curvature, and the Gauss map extends meromorphically to the conformal compactification.

\medskip

\section{Maximum number of tangency points of a horosphere packing}
\label{section-packing}
We start with the 2-dimensional case, namely horocycle packings in the hyperbolic
plane $\H^2$, because in this case we know the exact answer.
\begin{theorem}
For $n\geq 2$,
the maximum number of tangency points that a packing of $n$ horocycles may
have is $2n-3$.
\end{theorem}
Proof:
We work in the disk model of $\H^2$. 
Let us first see that the bound $2n-3$ can be achieved.
Start with two tangent horocycles.
Given two tangent horocycles, there exists two other horocycles which are tangent to
both of them, by Apolonius three circle theorem (with the ideal boundary of $\H^2$ as the
third circle). Choose one of them and iterate the process. The number of tangency
points is increased by 2 at each step, so this gives, for each $n\geq 2$,
a packing of $n$ horocycles with $2n-3$ tangency points.

\medskip

Next let us prove that one cannot do better. Consider an arbitrary packing of $n$
horocycles with $m$ tangency points.
We construct a planar graph with $n+1$ vertices and $m+n$ edges as follows:
The vertices are the $n$ points
$p_1,\cdots,p_n$ where the horocycles touch the unit circle (ideal boundary of
$\H^2$),
plus the point at $\infty$.
We connect two vertices $p_i$ and $p_j$ by the geodesic joining them if the corresponding
horocycles are tangent. 
We also connect each point $p_i$ to $\infty$ by a radial arc.
This gives a planar graph. Moreover, each of its faces has at least 3 edges in its boundary.
By a standard inequality for planar graph, the number of edges is at most 
$3N-6$ where $N$ is the number of vertices. This gives
$$m+n\leq 3(n+1)-6$$
which gives the result.
\cqfd
 
 \medskip
Let us now consider the 3-dimensional case of horosphere packings in $\H^3$.
Start with three mutually tangent horospheres.
Given three mutually tangent horospheres, there exists two other horospheres
which are tangent to all three (by the 3-dimensional version of Apolonius
theorem). Choose one of them and iterate the process. This gives, for each $n\geq 3$,
a packing of $n$ horospheres with $m=3n-6$ tangency points.
This is the construction described by Pacard and Pimentel in \cite{pacard}.
A packing obtained by this process is called an apolonian packing.

\medskip

Here is how to do better.
We work in the half-space model of $\H^3$.
Let $\Lambda$ be the equilateral lattice in the horizontal plane generated by the vectors $(1,0)$
and $(\cos\frac{\pi}{3},\sin\frac{\pi}{3})$. 
Consider some radius $R$. Let $\Lambda_R=\Lambda\cap\overline{D}(0,R)$.
For each point $p\in\Lambda_R$, consider
a euclidean sphere $S_p$ of radius $\frac{1}{2}$ centered at $(p,\frac{1}{2})$.
Stack on top of that the horizontal plane $x_3=1$, which is a horosphere in the half-space
model and is tangent to all horospheres $S_p$.
Let $n(R)$ be the total number of horospheres and $m(R)$ the number of tangency points.
For example, for $R=1$, we get $n(1)=8$ and $m(1)=19$, so this is already better 
than an apolonian packing.
For $R=2$, we get $n(2)=20$ and $m(2)=61$.
\medskip

Let us estimate the ratio $m(R)/n(R)$ for large values of $R$.
Each sphere $S_p$ is tangent to at most $6$ other spheres, so there are at most $3(n(R)-1)$ tangency points between them.
Therefore, adding the $n(R)-1$ tangency points with the horizontal horosphere,
$$m(R)\leq 4(n(R)-1).$$
On the other hand, if a sphere $S_p$ has less than
$6$ tangency points with the other spheres, then $p\in\Lambda_R\setminus D(0,R-1)$. The cardinal of this set is $O(R)$. Since $n(R)=O(R^2)$,
this gives
$$m(R)\geq 4n(R)-O(\sqrt{n(R)}).$$
Hence
$$\lim_{R\to\infty}\frac{m(R)}{n(R)}=4.$$
In general, we have the following bound.
\begin{theorem}
\label{theorem3}
For a packing of $n\geq 5$ horospheres, the number $m$ of tangency points
satisfies
$$m\leq 5n-16.$$
\end{theorem}
Proof. First consider a packing of $n=5$ horospheres. Then the five horospheres cannot
be pairwise tangent to each other.
Indeed, if this is the case, one can assume (by an isometry)
that one of the horosphere is the horizontal plane $x_3=1$. The statement then
boils down to the fact that four circles of the same radius in the plane cannot be
pairwise tangent to each other. This implies that $m\leq 9$ for $n=5$.
\medskip

Consider then a packing of $n>5$ horospheres. We work in the upper-half
space model of $\H^3$. We may choose the model so that all horospheres have
distinct euclidean radius. (Indeed, generically, an isometry will change the ratio of
the euclidean radii of two horospheres.) By the following lemma, the
horosphere which has the smallest euclidean radius has at most 5 tangency points.
Removing that horosphere, Theorem \ref{theorem3} follows by induction. 
\cqfd

\begin{lemma}
Consider three horospheres $S_1$, $S_2$, $S_3$ such that $S_1$, $S_2$ are tangent,
$S_1$, $S_3$ are tangent, and $S_2$, $S_3$ are either disjoint or tangent.
Let $R_i$ be the euclidean radius of $S_i$ and $(p_i,R_i)$ its center.
Assume that $R_1<R_2\leq R_3$. Then the angle $\theta_1$ of the triangle $p_1p_2p_3$
at $p_1$ satisfies $\theta_1>\frac{\pi}{3}$.
\end{lemma}
Proof:
If $S_i$ and $S_j$ are disjoint then the distance between their centers is greater than
$R_i+R_j$, hence
$$||p_i-p_j||^2+(R_i-R_j)^2>(R_i+R_j)^2$$
which gives 
$$||p_i-p_j||^2>4R_i R_j.$$
If $S_i$ and $S_j$ are tangent, these inequalities are equalities.
Trigonometry in the triangle $p_1p_2p_3$ gives
$$\cos\theta_1=\frac{||p_1-p_2||^2+||p_1-p_3||^2-||p_2-p_3||^2}
{2||p_1-p_2||\,||p_1-p_3||}\leq\frac{R_1 R_2+R_1 R_3 -R_2R_3}{2R_1\sqrt{R_2 R_3}}
%<\frac{R_1R_2}{2R_1\sqrt{R_2R_3}}
<\frac{1}{2}.$$
\cqfd

\medskip

Let $m_{\mbox{\scriptsize max}}(n)$ be the maximum number of tangency points that
a packing of $n$ horospheres may have.
Collecting the above results:
\begin{corollary}
$$\limsup_{n\to\infty} \frac{m_{\mbox{\scriptsize max}}(n)}{n}\in[4,5].$$
\end{corollary}

\section{Background material}
\label{section-preliminaries}
\subsection{Bryant representation}
\label{section-bryant}
\subsubsection{The Minkowski model of $\H^3$.}
\label{section-minkowski}
Let $\L^4$ denote the 4-dimensional lorentzian space, namely $\R^4$,
with coordinates denoted $x_0,x_1,x_2,x_3$ and metric $-dx_0^2+dx_1^2+dx_2^2+dx_3^2$. The hyperboloid
$$\{x\in\L^4\,:\, \langle x,x\rangle=-1,\; x_0>0\}$$
with the induced metric, is the Minkowski model of $\H^3$.
Bryant identifies $\L^4$ with the space of $2\times 2$ hermitian matrices
by identifying $(x_0,x_1,x_2,x_3)$ with the matrix
$$X=\left(\begin{array}{cc} x_0+x_3 & x_1+\i x_2\\
x_1-\i x_2 & x_0 -x_3\end{array}\right).$$
Then $\H^3$ becomes the set of positive definite hermitian matrices with 
determinant $1$.
\subsubsection{Bryant representation}
\label{section-bryant-representation}
Recall that a holomorphic map $F:\Sigma\to SL(2,\C)$ is null if $\det(F^{-1}dF)=0$.
\begin{theorem}[Bryant \cite{bryant}]
\label{theorem-bryant}
Let $\Sigma$ be a simply connected Riemann surface.
Let $F:\Sigma\to SL(2,\C)$ be a holomorphic
null immersion. Then $FF^*:\Sigma\to\H^3$
is a smooth conformal CMC-1 immersion.
Conversely, if  $f:\Sigma\to\H^3$ is a conformal CMC-1
immersion, there exists a holomorphic null immersion $F:\Sigma\to SL(2,\C)$
such that $f=FF^*$. Moreover, $F$ is unique up to right multiplication by
a constant matrix $H\in SU(2)$.
\end{theorem}
If $\Sigma$ is not simply connected, then $F$ is only well defined on the universal
cover of $\Sigma$ and will have $SU(2)$-valued monodromy.
\subsubsection{Global meromorphic data}
\label{section-global-data}
Assume we are given a CMC-$1$ immersion $f:\Sigma\to\H^3$.
By Theorem \ref{theorem-bryant}, we can write locally $f=FF^*$
where $F$ is a null holomorphic, $SL(2,\C)$-valued map.
Consider the matrix of holomorphic 1-forms
$$A(z)=(dF(z))F(z)^{-1}\in \mathfrak{sl}(2,\C)$$
where $\mathfrak{sl}(2,\C)$ is the Lie algebra of $2\times 2$ matrices
whose trace is zero.
Then $A$ is well defined on $\Sigma$: replacing $F$ by $FH$
does not change $A$.
Define $G=\frac{A_{11}}{A_{21}}$ and $\myomega=A_{21}$, where
$A_{ij}$ denote the coefficients of $A$. Then $G$ is a meromorphic function
and $\myomega$ is a holomorphic 1-form, both globally defined on $\Sigma$
(except in the exceptional case where $A_{21}\equiv 0$).
Since the trace and determinant of $A$ are zero, we have
\begin{equation}
\label{eq-A}
A=\left(\begin{array}{cc}
G & -G^2\\ 1 & -G\end{array}\right)\myomega
\end{equation}
The function $G$ is the Gauss map introduced
by Bryant, and $Q=\myomega\,dG$ is the Hopf quadratic differential.
We will see the geometric meaning of the Gauss map $G$ in the next section.
We call $(G,\myomega)$ the meromorphic data for the immersion $f$.

\medskip

Conversely, here is a recipe to construct $CMC$-1 immersions in $\H^3$.
Start with a Riemann surface $\Sigma$, a meromorphic function $G$ and a holomorphic 1-form $\myomega$ on $\Sigma$, such that $\myomega$ has a zero at each pole of $G$ with twice the multiplicity,
and has no other zeros. Define the matrix
$A$ by \eqref{eq-A}.
Then $A$ is a holomorphic, $\mathfrak{sl}(2,\C)$-valued 1-form on $\Sigma$
which does not vanish.
Solve the linear differential system
\begin{equation}
\label{eq-F}
dF(z)= A(z) F(z)
\end{equation}
with initial data $F(z_0)=F_0\in SL(2,\C)$.
The solution is a multi-valued holomorphic null immersion $F:\Sigma\to SL(2,\C)$.
(It is an immersion because $A(z)\neq 0$).
If its monodromy happens to be in $SU(2)$, then
the immersion $f=FF^*:\Sigma\to\H^3$ is well defined and has CMC-1.
\begin{remark}
Many authors consider instead the
matrix $F^{-1}dF$ and write
\begin{equation}
\label{eq-Weierstrass-local}
F^{-1}dF=\left(\begin{array}{cc} g & -g^2\\ 1& -g\end{array}\right)\omega
\end{equation}
where $g$ is a meromorphic function and $\omega$ is a holomorphic 1-form.
It turns out that $(g,\omega)$ is the Weierstrass data for the corresponding minimal
surface in $\R^3$ via the Lawson correspondence \cite{lawson}, which is why
\eqref{eq-Weierstrass-local} has become so popular.
The problem is that unless $\Sigma$ is simply connected,
$F^{-1}dF$ is not globally defined on $\Sigma$.
This is not surprising, since the Lawson correspondence is local.
So this point of view is not appropriate if we are to construct high genus examples.
\end{remark}
\subsubsection{Bryant representation in the half-space model}
\label{section-half-space}
A more familiar model of $\H^3$ is
the half-space $x_3>0$ in $\R^3$, with conformal metric
$x_3^{-2}(dx_1^2+dx_2^2+dx_3^2)$.
The following results are proved in \cite{rosenberg}.
An orientation preserving isometry from the Minkowski model to the half-space model is given by
$$\Phi(x_0,x_1,x_2,x_3)=\left(\frac{x_1}{x_0-x_3},\frac{x_2}{x_0-x_3},
\frac{1}{x_0-x_3}\right).$$
The immersion $\Phi\circ f$ is given in the half-space model by
\begin{equation}
\label{eq-half-space}
x_1+\i x_2=\frac{F_{11} \overline{F_{21}}+F_{12}\overline{F_{22}}}
{|F_{21}|^2+|F_{22}|^2}\qquad\qquad
x_3=\frac{1}{|F_{21}|^2+|F_{22}|^2}
\end{equation}
where $F_{ij}$ denote the coefficients of the matrix $F$.
The ideal boundary of $\H^3$ in this model is $\C\cup\{\infty\}$.
The Gauss map $G$ has the following geometric interpretation:
the normal geodesic ray originated from $\Phi\circ f(z)$ (in the direction of the mean curvature
vector) hits the ideal boundary at the point $G(z)$.
\subsubsection{Isometries}
\label{section-isometries}
The Lie group $SL(2,\C)$ acts isometrically on $\L^4$ by the representation
$$H\cdot X=HXH^*$$
where $H\in SL(2,\C)$ and $X\in \L^4$.
The action preserves $\H^3$ and its kernel is $\{\pm I_2\}$, so we recognize
$PSL(2,\C)$ as the group of direct isometries of $\H^3$.
The action of $SL(2,\C)$ on $\H^3$ extends to the ideal boundary as homographic
transformation of the Riemann sphere, namely
$$H\cdot z=\frac{H_{11} z+H_{12}}{H_{21}z+H_{22}}.$$
If $f:\Sigma\to\H^3$ is a conformal CMC-1 immersion with Gauss map $G$
and null holomorphic map $F$, 
then $H\cdot f$  has Gauss map $H\cdot G$ and null holomorphic map $HF$.
\subsubsection{Horospheres}
\label{section-horospheres}
The Gauss map is constant on a horosphere, and that constant is the limit point
of the horosphere on the ideal boundary of $\H^3$ in the half-space model.
(This follows from the geometric interpretation of the Gauss map in the half-space
model.)
If the horosphere is not a horizontal plane, then its meromorphic data is
$\Sigma=\C$,
 $G=c$,
 $\myomega=\lambda\,dz$
 where $c$ and $\lambda$ are complex constants.
 The constant $\lambda$ has no geometrical meaning and depends on the chosen
conformal parametrization of the horosphere.
The matrix $A$ is given by
$$A=\lambda\left(\begin{array}{cc} c & -c^2\\ 1 & -c\end{array}\right).$$
If the horosphere is a horizontal plane then $G\equiv \infty$
and we are in the exceptional case where the meromorphic data $(G,\myomega)$
is not defined.
In this case, one has
$$A=\left(\begin{array}{cc} 0 & \lambda\\ 0 & 0\end{array}\right)$$
for some complex number $\lambda$.
In any case, since the matrix $A$ is constant,
the solution to \eqref{eq-F} is $F(z)=\exp(zA)F_0$.
%\begin{remark}
%Let $p\in\C$ and $H=F(p)$. Then $HF(z)=F(z+p)$, so the matrix $H$
%represents an isometry of $\H^3$ which preserves the horosphere and
%maps $f(0)$ to $f(p)$.
%\end{remark}
\subsection{Linear differential systems}
\label{section-linear-systems}
In this section, we recall standard facts about linear differential systems on
a Riemann surface and setup some notations. A basic reference is
\cite{teschl}.
Let $\Sigma$ be a Riemann surface and $A$ an $n\times n$ matrix of
holomorphic 1-forms on $\Sigma$. We consider the first order linear differential system
on $\Sigma$
\begin{equation}
\label{eq-system}
dY(z)=A(z)Y(z)
\end{equation}
\subsubsection{Local theory: the principal solution}
\label{section-local-theory}
Assume that $\Sigma$ is simply connected.
Given $z_0\in \Sigma$, \eqref{eq-system} has a unique solution
$Y:\Sigma\to GL(n,\C)$ such that $Y(z_0)=I_n$.
Following \cite{teschl}, we write $Y(z)=\Pi(z,z_0)$.
The map $\Pi:\Sigma\times\Sigma\to GL(n,\C)$ is holomorphic in both variables and is
called the principal solution. It satisfies
$$\Pi(z_3,z_2)\Pi(z_2,z_1)=\Pi(z_3,z_1).$$
Given $Y_0\in GL(n,\C)$, the solution $Y$ such that $Y(z_0)=Y_0$ is given by
$Y(z)=\Pi(z,z_0)Y_0$.
\subsubsection{Global theory}
\label{section-global-theory}
Now assume that $\Sigma$ is not simply connected. Then the principal solution
$\Pi(z,z_0)$ is not well defined: it depends on the homotopy class of the path from
$z_0$ to $z$. If $\gamma:[0,1]\to\Sigma$ is a path from $z_0$ to $z$,
the solution $Y$
of \eqref{eq-system} such that $Y(z_0)=I_n$, which exists in a simply connected
neighborhood of $\gamma(0)$, can be analytically continued along $\gamma$.
Its value at $\gamma(1)$ will be denoted $\Pi(\gamma)$. When the path $\gamma$
is clear from the context, we will still use the notation $\Pi(z,z_0)$.
\medskip

Let $\gamma_1,\gamma_2:[0,1]\to\Sigma$ be two paths such that $\gamma_1(1)=\gamma_2(0)$.
We denote by $\gamma_2\cdot\gamma_1$ the path obtained by composing 
$\gamma_1$ and $\gamma_2$.
(The usual notation is $\gamma_1\cdot\gamma_2$ but in this context, it is more
convenient, and customary, to reverse the order).
Then
\begin{equation}
\label{eq-morphism}\Pi(\gamma_2\cdot\gamma_1)=\Pi(\gamma_2)\Pi(\gamma_1).
\end{equation}
In particular, $\Pi$ is a morphism from the fundamental group
$\pi_1(\Sigma,z_0)$ to $GL(n,\C)$.
(With the usual notation for the product in the fundamental group, it would be an anti-morphism.)
\subsubsection{Monodromy}
\label{section-intro-monodromy}
The monodromy of a solution is usually defined as follows. Let $Y$ be a solution
of \eqref{eq-system} in a simply connected neighorhood $U$ of $z_0$.
Let $\gamma\in\pi_1(\Sigma,z_0)$. Analytic continuation of $Y$ along $\gamma$
gives another solution of \eqref{eq-system} in $U$, which is denoted
$\gamma\cdot Y$. There exists a matrix $M\in GL(n,\C)$ such that $\gamma\cdot Y=YM$.
The matrix $M$ is called the monodromy of $Y$ along $\gamma$ and is denoted
$M_{\gamma}(Y)$. In term of the principal solution, one has
\begin{equation}
\label{eq-monodromy-conjugate}
M_{\gamma}(Y)=Y(z_0)^{-1} \Pi(\gamma) Y(z_0).
\end{equation}
\subsection{Opening nodes}
\label{section-opening-nodes}
We recall the standard construction of opening nodes.
Consider $n$ copies of the Riemann sphere $\CC=\C\cup\{\infty\}$, labelled
$\CC_1,\cdots,\CC_n$. Consider $2m$ distinct points $p_1,\cdots,p_m$,
$q_1,\cdots,q_m$ in the disjoint union $\CC_1\cup\cdots\cup\CC_n$.
Identify $p_i$ with $q_i$ for $1\leq i\leq m$. This defines a Rieman surface with
nodes which we denote $\Sigma_0$. We assume $\Sigma_0$ is connected.
\medskip

To open nodes, consider local complex coordinates $v_i:V_i\to D(0,1)$ in a neighborhood
of $p_i$ and $w_i:W_i\to D(0,1)$ in a neighborhood of $q_i$, with
$v_i(p_i)=0$ and $w_i(q_i)=0$. We assume that the neighborhoods
$V_1,\cdots,V_m,W_1,\cdots,W_m$ are disjoint in $\CC_1\cup\cdots\cup\CC_n$.
Consider, for each $1\leq i\leq m$, a complexe parameter
$t_i$ with $|t_i|< 1$.
If $t_i=0$, identify $p_i$ with $q_i$ as above.
If $t_i\neq 0$, remove the disks $|v_i|\leq |t_i|$ and $|w_i|\leq |t_i|$. Identify
the point $z\in V_i$ with the point $z'\in W_i$ such that
$$v_i(z)w_i(z')=t_i.$$
This creates a Riemann surface, possibly with nodes,
which we denote $\Sigma_{\bft}$, where $\bft=(t_1,\cdots,t_n)$.
When all $t_i$ are non zero, $\Sigma_{\bft}$ is a genuine Riemann surface of
genus $g=m-n+1$.

\medskip

Observe that if $t_i\neq 0$, the circle $|v_i|=1$ is homologous in $\Sigma$ to
the circle $|w_i|=1$ with the opposite orientation. Consequently,
if $\omega$ is a holomorphic 1-form on $\Sigma_{\bf t}$,
\begin{equation}
\label{eq-opposite}
\int_{|v_i|=1}\omega=-\int_{|w_i|=1}\omega.
\end{equation}
This makes the following definition natural.
\begin{definition}[Bers]
A regular differential on $\Sigma_{\bft}$ is a holomorphic 1-form, which is allowed
to have simples poles at $p_i$ and $q_i$ if $t_i=0$, with opposite residues.
\end{definition}
By a theorem of Fay \cite{fay}, the space of regular differentials on $\Sigma_{\bft}$
has dimension $g$ and admits a basis which depends holomorphically on $\bft$
in a neighborhood of $0$.
For $1\leq i\leq n$, let $J_i^+$ and $J_i^-$ be the set of indices $j$ such that
$p_j\in\CC_i$ and $q_j\in\CC_i$, respectively.
As a consequence of Fay's theorem, one has:
\begin{theorem}
\label{theorem-opening-nodes}
For $\bft$ in a neighborhood of $0$, and for $\bfa=(a_j)_{1\leq j\leq m}\in\C^m$ satisfying
\begin{equation}
\label{eq-homology}
\sum_{j\in J_i^+}a_j-\sum_{j\in J_i^-}a_j=0\qquad\mbox{ for $1\leq i\leq n$}
\end{equation}
there exists a unique regular differential $\omega=\omega_{\bft,\bfa}$ on $\Sigma_{\bft}$
such that
\begin{equation}
\label{eq-prescribe-periods}
\int_{|v_j|=1}\omega=a_j\qquad \mbox{ for $1\leq j\leq m$}.
\end{equation}
Moreover, $\omega_{\bft,\bfa}$ depends holomorphically on ${\bft}$ (away from the nodes).
\end{theorem}
Proof: from Cauchy theorem in $\CC_i$ and \eqref{eq-opposite}, we see
that \eqref{eq-homology} is necessary for $\omega$ to exist.
If $\bft=0$, the map $\omega\mapsto(\int_{|v_j|=1}\omega)_{1\leq j\leq m}$ is an isomorphism
from the space of regular differentials on $\Sigma_0$ to the space of vectors
$\bfa\in\C^m$ satisfying \eqref{eq-homology}. (This follows from the fact that
a holomorphic 1-form on the Riemann sphere with simple poles is entirely
determined by its residues.)
Using Fay's theorem, this remains true for $\bft$ in a neighborhood of $0$.
\cqfd

\medskip

Fay's proof is rather abstract and non-constructive. For an elementary proof of
Theorem \ref{theorem-opening-nodes} based on the contraction mapping principle,
see \cite{crelle}.
One has a similar result for meromorphic 1-forms with poles at some points
$r_1,\cdots,r_k$ distinct from the nodes:
\begin{theorem}
\label{theorem-meromorphic}
For $\bft$ in a neighborhood of $0$, one can define a regular meromorphic differential
$\omega$ on $\Sigma_{\bft}$
with poles at the points $r_1,\cdots,r_k$ by prescribing its
principal part at each pole and its periods as in \eqref{eq-prescribe-periods},
replacing \eqref{eq-homology} by the restriction coming from the residue theorem, namely
\begin{equation}
\label{eq-residue}
\sum_{j\in J_i^+}a_j-\sum_{j\in J_i^-}a_j
+2\pi \i\sum_{r_j\in \CC_i}\Res_{r_j} \omega=0\qquad\mbox{ for $1\leq i\leq n$}.
\end{equation}
\end{theorem}
Recall that the principal part of a meromorphic 1-form $\omega$ at a pole $r$
is its equivalence class under the relation: $\omega\sim\omega'$ if
$\omega-\omega'$ is holomorphic in a neighborhood of $r$.
The analogue of Fay's theorem for meromorphic differentials with simple poles is
proved by Masur in \cite{masur}. For a proof of Theorem \ref{theorem-meromorphic} in the case of poles of arbitrary order, see \cite{crelle}.
We will also need the following result to compute the partial derivatives of $\omega$ with
respect to $\bft$.
\begin{theorem}
\label{theorem-derivee-t}
The partial derivative $\frac{\partial}{\partial t_i}\omega_{\bft,\bfa}$ at $\bft=0$
has two double poles at $p_i$ and $q_i$, with principal parts
$$\frac{-dv_i}{v_i^2}\Res_{q_i}\frac{\omega_{0,\bfa}}{w_i}\qquad
\mbox{ at $p_i$}$$
$$\frac{-dw_i}{w_i^2}\Res_{p_i}\frac{\omega_{0,\bfa}}{v_i}\qquad
\mbox{ at $q_i$}$$
and has vanishing periods on all circles $|v_j|=1$.
\end{theorem}
This is proved in \cite{triply}, Lemma 3. See also \cite{crelle}, Remark 5.6.

\section{The meromorphic data $(\Sigma,G,\myomega)$}
\label{section-data}
\subsection{Notations}
The horospheres of our given horosphere packing are denoted $S_1,\cdots,S_n$.
We define the following sets, which we use to index various quantities:
$$I=\{(i,j)\in\[1,n\]^2\;:\;i<j \mbox{ and $S_i,S_j$ are tangent}\}$$
$$J_i^+=\{j\,:\, (i,j)\in I\}\qquad J_i^-=\{j\;:\;(j,i)\in I\}\qquad J_i=J_i^+\cup J_i^-.$$
Without loss of generality, we may assume (by applying an isometry) that each horosphere $S_i$ is not
a horizontal plane in the half-space model.
We fix a conformal parametrisation $f_i:\C\to S_i$ and let
$$G=c_i\qquad\myomega=\lambda_i dz$$
be its meromorphic data (see Section \ref{section-horospheres}).
For $j\in J_i$, we denote $p_{ij}$ the point in $\C$ such that
$f_i(p_{ij})$ is the point $S_i\cap S_j$.
Without loss of generality, we may assume (by changing the parametrization $f_i)$
that for each $i$, the disks
$D(p_{ij},1)$ for $j\in J_i$ and $D(0,1)$ are pairwise disjoint.
\subsection{The Riemann surface $\Sigma$}
We consider $n$ copies of the complex plane, denoted $\C_1,\cdots,\C_n$,
and $m$ copies of the Riemann sphere $\C\cup\{\infty\}$, denoted
$\CC_{ij}$ for $(i,j)\in I$.
We think of $p_{ij}$ as a point in $\C_i$.
%The points $0$ and $\infty$ in $\C_i$ will be denoted respectively $0_i$ and $\infty_i$.
The points $0$, $1$ and $\infty$ in $\CC_{ij}$ will be
denoted respectively $0_{ij}$, $1_{ij}$ and $\infty_{ij}$.
Heuristically, the reader should think of $\C_1,\cdots,\C_n$ as
the parametrization domain for
the horospheres $S_1,\cdots,S_n$, and $\CC_{ij}\setminus\{0,\infty\}$ for $(i,j)\in I$ as the parametrization domain for the
catenoidal necks connecting them.

\medskip
We define a Riemann surface  $\Sigma_0$ with $2m$ nodes by identifying $p_{ij}$ with
$0_{ij}$ and $p_{ji}$ with $\infty_{ij}$, for $(i,j)\in I$ and call it $\Sigma_0$.
To open nodes, we consider the natural coordinates $v_{ij}=z-p_{ij}$ in a neighborhood of $p_{ij}$ in $\C_i$,
$w_{ij}=z$ in a neighborhood of $0_{ij}$ in $\CC_{ij}$,
$v_{ji}=z-p_{ji}$ in a neighborhood of $p_{ji}$ in $\C_j$ and
$w_{ji}=\frac{1}{z}$ in a neighborhood of $\infty_{ij}$ in $\CC_{ij}$.
We open nodes as explained in Section \ref{section-opening-nodes},
introducing a complex parameter $t_{ij}$ to open the node
$p_{ij}\sim 0_{ij}$ and another parameter $t_{ji}$ to open the node $p_{ji}\sim \infty_{ij}$.
Let $\bft=(t_{ij},t_{ji})_{(i,j)\in I}$ be the collection of these parameters.
This defines a Riemann surface (possibly with nodes) which we denote $\Sigma$, or
$\Sigma_{\bft}$ when we need to emphasize the
dependance on the parameter $\bft$.
We denote $\overline{\Sigma}$ the compactification of $\Sigma$ obtained by adding the points $\infty_1,\cdots,\infty_n$,
where $\infty_i$ denotes the point at infinity in $\C_i$.
\medskip

For $(i,j)\in I$, we define $\alpha_{ij}$ as the homology class of the circle
$|z-p_{ij}|=1$ in $\C_i$, with the positive orientation. This is homologous, in
$\Sigma$, to the unit circle in $\CC_{ij}$, with the negative orientation, and
to the circle $|z-p_{ji}|=1$ in $\C_j$, also with the negative orientation.
\subsection{The Gauss map $G$}
Here are our requirements on the Gauss map $G$. At $\infty_i$, it should take the value
$c_i$. It should have a simple pole in
each Riemann sphere $\CC_{ij}$ (because on each catenoidal neck, we expect
a point where the mean curvature vector is vertical pointing up, in the half-space
model). We choose the identification of the Riemann sphere with $\C\cup\{\infty\}$
so that this pole is $z=1$.
The following proposition tells us that these requirements completely determine $G$.
\begin{proposition}
\label{proposition-gauss}
For $\bft$ small enough, there exists a unique meromorphic
function $G=G_{\bft}$ on $\overline{\Sigma}_{\bft}$ with the following properties :
\begin{itemize}
\item $G_{\bft}$ has $m$ simple poles at the points $1_{ij}$ for $(i,j)\in I$,
\item $G_{\bft}(\infty_i)=c_i$ for $i=1,\cdots,n$.
\end{itemize}
Moreover, $G_{\bft}$ depends holomorphically on $\bft$
(away from the nodes and its poles) and at $\bft=0$, we have
\begin{equation}
\label{eq-G0}
G_0(z)=\left\{\begin{array}{ll}
c_i & \mbox{ in $\C_i$}\\
c_j+\displaystyle\frac{c_j-c_i}{z-1} & \mbox{ in $\CC_{ij}$}
\end{array}\right.
\end{equation}
\end{proposition}

Proof: We first define the differential $\mu=dG$ and we recover $G$ by integration.
By Theorem \ref{theorem-meromorphic}, there exists a unique meromorphic differential $\mu_{\bft}$ on $\overline{\Sigma}_{\bft}$
which has $m$ double poles at $1_{ij}$ for $(i,j)\in I$
with principal part
$$r_{ij}\frac{dz}{(z-1)^2}$$
and has vanishing period on all cycles $\alpha_{ij}$.
Here the $m$ complex numbers $r_{ij}$ are free parameters.
At $\bft=0$, we have
$$\mu_0=\left\{\begin{array}{ll}
0 & \mbox{ in $\C_i$}\\
r_{ij}\frac{dz}{(z-1)^2} &\mbox{ in $\CC_{ij}$.}
\end{array}\right.$$
\begin{lemma}
\label{lemma-gauss}
For $\bft$ in a neighborhood of $0$, there exist unique values of the parameters
$r_{ij}$ such that
\begin{equation}
\label{eq-integral-mu}\int_{\infty_i}^{\infty_j}\mu_{\bft}=c_j-c_i
\qquad  \mbox{ for $(i,j)\in I$}.
\end{equation}
Moreover, each $r_{ij}$ depends holomorphically on $\bft$, and when $\bft=0$,
$r_{ij}=c_i-c_j$.
\end{lemma}
Proof: First observe that \eqref{eq-integral-mu} is a linear system of $m$ linear equations
with $m$ unknowns $r_{ij}$, $(i,j)\in I$. Hence it suffices to prove that
the system \eqref{eq-integral-mu} is invertible when $\bft=0$.
In \eqref{eq-integral-mu}, it is understood that the path from $\infty_i$
to $\infty_j$ goes through $\CC_{ij}$. There is no canonical way to choose this
path, but all choices are homologous modulo $\alpha_{ij}$.
Since
$\mu$ has no period on $\alpha_{ij}$, $\int_{\infty_i}^{\infty_j}\mu_{\bft}$
is a well defined holomorphic function of $\bft$. Moreover, by Lemma 4 in
\cite{triply}, this function extends holomorphically at $\bft=0$ with value
$$\int_{\infty_i}^{\infty_j}\mu_0=
\int_{\infty_i}^{p_{ij}} 0+\int_{0_{ij}}^{\infty_{ij}}r_{ij}\frac{dz}{(z-1)^2}+
\int_{p_{ji}}^{\infty_j}0
=-r_{ij}.$$
The Lemma follows.
\cqfd

\medskip

Returning to the proof of Proposition \ref{proposition-gauss}, we define the function
$G_{\bft}$ on $\overline{\Sigma}_{\bft}$ by
$$G_{\bft}(z)=c_1+\int_{\infty_1}^z\mu_{\bft}.$$
By Lemma \ref{lemma-gauss}, and the fact that $\mu_{\bft}$ has no residues
and no period on the cycles $\alpha_{ij}$, $G_{\bft}$ is well defined on $\overline{\Sigma}_{\bft}$
(meaning that the integral does not depend on the path from $\infty_1$ to $z$)
and has all desired properties.
\cqfd
\subsection{The holomorphic differential $\myomega$}
Here are our requirements on the holomorphic differential $\Omega$. It should
have a double pole at $\infty$, with leading term $\lambda_i dz$, just like the
horosphere $S_i$. It also needs a double zero at each pole of $G$.
We define $\Omega$ by prescribing poles, principal parts and periods,
using Theorem \ref{theorem-meromorphic}. Then we adjust the parameter
$\bft$ so that $\Omega$ has the required zeros.
\begin{definition}
Consider $m$ complex parameters $a_{ij}$, $(i,j)\in I$ and let
$\bfa=(a_{ij})_{(i,j)\in I}$.
We define $\myomega=\myomega_{\bft,\bfa}$ as the unique meromorphic 1-form on $\overline{\Sigma}_{\bft}$ with $n$
double poles at $\infty_1,\cdots,\infty_n$, with principal part
$$\lambda_i dz+\sum_{j\in J_i^+}a_{ij}\frac{dz}{z}-\sum_{j\in J_i^-}a_{ji}\frac{dz}{z} \qquad\mbox{ at $\infty_i$}$$
and periods
$$\qquad \int_{\alpha_{ij}}\myomega=2\pi \i\,a_{ij}
\qquad \mbox{ for $(i,j)\in I$}.$$
It depends holomorphically (away from its poles and the nodes) on $\bft$.
Moreover, at $\bft=0$ we have
\begin{equation}
\label{eq-omega0}
\myomega_{0,\bfa}=\left\{\begin{array}{ll}
\displaystyle\lambda_i dz+\sum_{j\in J_i^+}a_{ij}\frac{dz}{z-p_{ij}}
-\sum_{j\in J_i^-}a_{ji}\frac{dz}{z-p_{ji}} &\mbox{ in $\C_i$}\\
\displaystyle -a_{ij}\frac{dz}{z}&\mbox{ in $\CC_{ij}$}\end{array}\right.
\end{equation}
\end{definition}
Note that the residue of the prescribed principal part at $\infty_i$ is
$$\Res_{\infty_i}\myomega=-\sum_{j\in J_i^+}a_{ij}+\sum_{j\in J_i^-}a_{ji}$$
so Equation \eqref{eq-residue} holds.
\begin{proposition}
\label{proposition-zeros}
For $\bfa$ in a neighborhood of $0$, there exists a unique value $\bft(\bfa)$, depending
holomorphically on $\bfa$, such that $\myomega_{\bft(\bfa),\bfa}$ has a double zero at
each pole of $G$, and has no other zeros in $\Sigma$ (provided all parameters
$a_{ij}$ are non-zero).
Moreover, for each $(i,j)\in I$, we have
\begin{equation}
\label{eq-aij-nul}
a_{ij}=0\quad\Rightarrow\quad t_{ij}(\bfa)=t_{ji}(\bfa)=0
\end{equation}
\begin{equation}
\label{eq-derivee-tij}
\frac{\partial t_{ij}(\bfa)}{\partial a_{ij}}|_{\bfa=0}=\frac{-1}{2\lambda_i}
\end{equation}
\begin{equation}
\label{eq-derivee-tji}
\frac{\partial t_{ji}(\bfa)}{\partial a_{ij}}|_{\bfa=0}=\frac{1}{2\lambda_j}.
\end{equation}
\end{proposition}
Proof: let $(i,j)\in I$.
Using Theorem \ref{theorem-derivee-t} and \eqref{eq-omega0}, we compute the partial derivatives
of $\myomega_{\bft,\bfa}$ in $\CC_{ij}$ at $(\bft,\bfa)=(0,0)$.
\begin{eqnarray}
\label{eq-derivee-omega1}
\frac{\partial\myomega_{\bft,\bfa}}{\partial t_{ij}}&=&
\frac{-dz}{z^2}\Res_{p_{ij}}\frac{\myomega_{0,0}}{z-p_{ij}}=-\lambda_i\frac{dz}{z^2}
\\
\label{eq-derivee-omega2}
\frac{\partial\myomega_{\bft,\bfa}}{\partial t_{ji}}&=&
\frac{-dw_{ji}}{w_{ji}^2}\Res_{p_{ji}}\frac{\myomega_{0,0}}{z-p_{ji}}=\lambda_j\,dz
\qquad\mbox{ (recall $w_{ji}=\displaystyle\frac{1}{z}$)}\\
\label{eq-derivee-omega3}
\frac{\partial\myomega_{\bft,\bfa}}{\partial a_{ij}}&=&
\frac{-dz}{z}
\end{eqnarray}
The partial derivatives of $\myomega_{\bft,\bfa}$ in $\CC_{ij}$ with respect to all other
parameters $t_{k\ell}$ and $a_{k\ell}$ are zero.
Write $\myomega_{\bft,\bfa}=f_{ij}(\bft,\bfa,z) dz$ in $\CC_{ij}$. We want to solve
$f_{ij}(\bft,\bfa,1)=f'_{ij}(\bft,\bfa,1)=0$.
The Jacobian matrix of $(f_{ij}(\bft,\bfa,1),f'_{ij}(\bft,\bfa,1))$ with respect to $(t_{ij},t_{ji})$
is
$$\left(\begin{array}{cc}
-\lambda_i & \lambda_j\\
2\lambda_i & 0\end{array}\right).$$
The existence of the solution $\bft(\bfa)$ then follows from the implicit function theorem
applied to the map $(f_{ij}(\bft,\bfa,1),f'_{ij}(\bft,\bfa,1))_{(i,j)\in I}$ whose Jacobian has block diagonal form.
Next consider some $(i,j)\in I$ and assume that $a_{ij}=0$. If $t_{ij}=t_{ji}=0$, then $\Omega_{\bft,\bfa}$ has two
simple poles at $0_{ij}$ and $\infty_{ij}$, with residue $\pm a_{ij}=0$, hence is
holomorphic in $\CC_{ij}$, so $\Omega_{\bft,\bfa}\equiv 0$ in $\CC_{ij}$. So \eqref{eq-aij-nul}
follows from uniqueness in the implicit function theorem.
Equations \eqref{eq-derivee-tij} and \eqref{eq-derivee-tji} are obtained by differentiating
$f_{ij}(\bft(\bfa),\bfa,1)=f_{ij}'(\bft(\bfa),\bfa)=0$ with respect to $a_{ij}$, using
\eqref{eq-derivee-omega1}, \eqref{eq-derivee-omega2} and \eqref{eq-derivee-omega3}.
Finally, if all parameters $a_{ij}$ are non-zero, then
\eqref{eq-derivee-tij} and \eqref{eq-derivee-tji} imply that all parameters
$t_{ij}$ and $t_{ji}$ are non-zero, so $\Sigma$ is a genuine compact Riemann surface
of genus $g=m-n+1$.
$\myomega$ is a meromorphic 1-form with $n$ doubles poles, so its number of zeros,
counting multiplicity, is $2n+2g-2=2m$. This ensures that $\myomega$ has no
other zeros than the double zeros it has at the $m$ poles of $G$.
\cqfd
\subsection{Partial derivatives with respect to the parameter $a_{ij}$}
\begin{proposition} For $(i,j)\in I$ we have at $\bfa=0$
\begin{equation}
\label{eq-derivee-G}
\frac{\partial G_{\bft(\bfa)}}{\partial a_{ij}}=\left\{\begin{array}{ll}
\displaystyle\frac{c_j-c_i}{2\lambda_i}\frac{1}{z-p_{ij}} &\mbox{ in $\C_i$}\\
\displaystyle\frac{c_j-c_i}{2\lambda_j}\frac{1}{z-p_{ji}} & \mbox{ in $\C_j$}\\
0 & \mbox{ in $\C_k$ for $k\neq i,j$}.\end{array}\right.
\end{equation}
\begin{equation}
\label{eq-derivee-omega}
\frac{\partial\myomega_{\bft(\bfa),\bfa}}{\partial a_{ij}}=\left\{\begin{array}{ll}
\displaystyle\frac{dz}{z-p_{ij}} & \mbox{ in $\C_i$}\\
\displaystyle\frac{-dz}{z-p_{ji}} & \mbox{ in $\C_j$}\\
\displaystyle\frac{(1-z)^2}{2z^2}dz & \mbox{ in $\CC_{ij}$}\\
0 & \mbox{ in $\C_k$ for $k\neq i,j$ and in $\CC_{k\ell}$ for $(k,\ell)\neq (i,j)$}.
\end{array}\right.
\end{equation}
\end{proposition}
Proof: Recall that $\mu=dG$. By Theorem \ref{theorem-derivee-t},
$$\frac{\partial \mu_{\bft}}{\partial t_{ij}}
=-\frac{dz}{(z-p_{ij})^2}\Res_{0_{ij}}\mu_0
=-(c_i-c_j)\frac{dz}{(z-p_{ij})^2}\qquad\mbox{ in $\C_i$.}$$
Hence, by the chain rule and \eqref{eq-derivee-tij},
\begin{equation}
\label{eq-derivee-dG}
\frac{\partial\mu_{\bft(\bfa)}}{\partial a_{ij}}=\frac{c_i-c_j}{2\lambda_i}\frac{dz}{(z-p_{ij})^2}
\qquad\mbox{ in $\C_i$}.
\end{equation}
Integrating, we obtain the first line of \eqref{eq-derivee-G}. The proof of the second line is
entirely similar. Regarding \eqref{eq-derivee-omega}, we have (using again Theorem \ref{theorem-derivee-t})
$$\frac{\partial\myomega_{\bft,\bfa}}{\partial t_{ij}}=0
\qquad\qquad
\frac{\partial\myomega_{\bft,\bfa}}{\partial a_{ij}}=\frac{dz}{z-p_{ij}}
\qquad\mbox{ in $\C_i$}.$$
The first line of \eqref{eq-derivee-omega} follows from the chain rule.
The proof of the second line is similar.
Using Equations
\eqref{eq-derivee-omega1}, \eqref{eq-derivee-omega2}, \eqref{eq-derivee-omega3}
and the chain rule, we have
$$\frac{\partial\myomega_{\bft(\bfa),\bfa}}{\partial a_{ij}}
=-\lambda_i\frac{dz}{z^2}\times\left(\frac{-1}{2\lambda_i}\right)
+\lambda_j\,dz\times\left(\frac{1}{2\lambda_j}\right)-\frac{dz}{z}
=\frac{(1-z)^2}{2z^2}dz\qquad\mbox{ in $\CC_{ij}$}.$$
\cqfd
\section{The monodromy problem}
\label{section-monodromy}
\subsection{Formulation of the problem}
\label{section-formulation}
We consider the matrix $A=A_{\bfa}$ defined by \eqref{eq-A}
with
$G=G_{\bft(\bfa)}$ and $\myomega=\myomega_{\bft(\bfa),\bfa}$.
Each coefficient of $A$ is a holomorphic differential on $\Sigma=\Sigma_{\bft(\bfa)}$.
Let $0_i$ be the point $z=0$ in $\C_i$.
Let $F:\Sigma\to SL(2,\C)$ be the solution of $dF=AF$ with initial condition
$F(0_1)=\Lambda_1$, where $\Lambda_1\in SL(2,\C)$ is a matrix we can
prescribe. (Observe that $F(z)\in SL(2,\C)$ because $A(z)\in\mathfrak{sl}(2,\C)$.)
The solution $F$ is of course only well defined on the universal cover of
$\Sigma$. We need to adjust the parameters
so that $F$ has $SU(2)$-valued monodromy, so $f=F F^*$ is well defined
on $\Sigma$.
Taking $0_1$ as a base point for the fundamental group and using
\eqref{eq-monodromy-conjugate}, this
is equivalent to
\begin{equation}
\label{eq-monodromy}
\forall\gamma\in\pi_1(\Sigma,0_1),\qquad \Lambda_1^{-1}\Pi(\gamma)\Lambda_1\in SU(2)
\end{equation}
where $\Pi$ denotes the principal solution of $dF=AF$ on $\Sigma$
(see Section \ref{section-local-theory}).
\medskip

Instead of using a set of generators of $\pi_1(\Sigma,0_1)$, which would involve in
a complicated way the ``combinatorics'' of our given horosphere packing, we
reformulate the monodromy problem in a more ``local'' way as follows.
For $(i,j)\in I$, let $\gamma_{ij}\in \pi_1(\Sigma,0_i)$ be a loop in $\C_i$ with base point $0_i$ which
goes around $p_{ij}$ and does not encircle any other node. We also define $\Gamma_{ji}$ as a path connecting $0_i$
to $0_j$ through $\CC_{ij}$ (to be defined more precisely later on).
\begin{proposition}
\label{proposition-formulation}
Given $n$ matrices $\Lambda_1,\cdots,\Lambda_n$ in $SL(2,\C)$, assume that
for all $(i,j)\in I$:
\begin{equation}
\label{eq-monodromy1}
\Lambda_i^{-1} \Pi(\gamma_{ij}) \Lambda_i\in SU(2)
\end{equation}
\begin{equation}
\label{eq-monodromy2}
\Lambda_j^{-1} \Pi(\Gamma_{ji})\Lambda_i\in SU(2)
\end{equation}
Then \eqref{eq-monodromy} is satisfied.
Moreover, the solution $F$ of $dF=AF$ with initial condition $F(0_1)=\Lambda_1$
satisfies $F(0_i)\in \Lambda_i\times SU(2)$, hence
$f(0_i)=\Lambda_i\Lambda_i^*$ for $1\leq i\leq n$.
\end{proposition}
Proof: Assume that \eqref{eq-monodromy1} and \eqref{eq-monodromy2}
hold for all $(i,j)\in I$.
Let $\Gamma_{ij}=\Gamma_{ji}^{-1}$. Using \eqref{eq-morphism}
$$\Lambda_i^{-1}\Pi(\Gamma_{ij})\Lambda_j=\left(
\Lambda_j^{-1}\Pi(\Gamma_{ji})\Lambda_i\right)^{-1}\in SU(2).$$
Define $\gamma_{ji}=\Gamma_{ji}\gamma_{ij}\Gamma_{ij}\in\pi_1(\Sigma,0_j)$.
Using \eqref{eq-morphism} again,
$$\Lambda_j^{-1}\Pi(\gamma_{ji})\Lambda_j=
\left(\Lambda_j^{-1}\Pi(\Gamma_{ji})\Lambda_i\right)
\left(\Lambda_i^{-1}\Pi(\gamma_{ij})\Lambda_i\right)
\left(\Lambda_i^{-1}\Pi(\Gamma_{ij})\Lambda_j\right)\in SU(2).$$
In other words, \eqref{eq-monodromy1} and \eqref{eq-monodromy2} also hold for
$(j,i)\in I$.
Let $\C_i^{\star}$ be $\C_i$ minus the disks $D(p_{ij},1)$ for $j\in J_i$.
The fundamental group $\pi_1(\C_i^{\star},0_i)$ is the free group with generators
$\gamma_{ij}$ for $j\in J_i$. Hence \eqref{eq-monodromy1} implies that
\begin{equation}
\label{eq-monodromy3}
\forall \delta\in\pi_1(\C_i^{\star},0_i)\quad \Lambda_i^{-1}\Pi(\delta)\Lambda_i
\in SU(2).
\end{equation}
Any element $\gamma\in\pi_1(\Sigma,0_1)$ is homotopic to a product of the
form
$$\delta_k\Gamma_{i_k\,i_{k-1}}\delta_{k-1}\Gamma_{i_{k-1}\,i_{k-2}}\delta_{k-2}\cdots
\delta_2\Gamma_{i_2\,i_1}\delta_1$$
where $k\in\N^*$, $i_1=i_{k}=1$ and $\delta_j\in\pi_1(\C_{i_j}^{\star},0_{i_j})$
for $1\leq j\leq k$.
So \eqref{eq-monodromy} follows from \eqref{eq-monodromy2} and \eqref{eq-monodromy3}.
\cqfd
\subsection{Choice of the matrices $\Lambda_1,\cdots\Lambda_n$.}
\label{section-choice-Lambda}
By the last statement of Proposition \ref{proposition-formulation},
choosing the matrices $\Lambda_i$ amounts
to prescribe the image of the points $0_1,\cdots,0_n$. Recall from the introduction
that we want to ``deflate'' the horosphere $S_i$ at speed $\xi_i$. This suggests the
following choice.
Consider $n$ fixed, positive numbers $\xi_1,\cdots,\xi_n$ and a real parameter $s$.
(These are the same parameters as in the introduction).
Let $O_i=f_i(0)\in S_i$, where
$f_i$ is the chosen conformal parametrization of the horosphere $S_i$.
Choose $\Lambda_i(s)\in SL(2,\C)$ so that 
$s\in[0,\infty)\mapsto\Lambda_i(s)\Lambda_i(s)^*$ (in the Minkowski model)
is the parametrization at speed $\xi_i$
of the geodesic ray normal to the horosphere $S_i$ at the point $O_i$ (in the direction
of the mean curvature vector).
The matrix $\Lambda_i(s)$ is unique up to right multiplication by an element in $SU(2)$, which is clearly irrelevant for the monodromy problem.
\subsection{Main result}
\label{section-main-result}
To solve the monodromy problem, we need to adjust the complex parameters $p_{ij}$ and
$p_{ji}$ for
$(i,j)\in I$. We will denote $p_{ij}^0$ and $p_{ji}^0$ the value of these parameters corresponding to the given horosphere packing
(namely, such that $f_i(p_{ij}^0)=f_j(p_{ji}^0)=S_i\cap S_j$).
The matrix of holomorphic 1-forms $A$ depends holomorphically on the parameters
$\bfa=(a_{ij})_{(i,j)\in I}$ and $\bfp=(p_{ij},p_{ji})_{(i,j)\in I}$ and will be denoted $A_{\bfa,\bfp}$.
The principal solution of $dF=A_{\bfa,\bfp} F$ will be denoted $\Pi_{\bfa,\bfp}$.
Our goal is to prove:
\begin{proposition}[solution of the monodromy problem]
\label{proposition-monodromy}
For $s>0$ small enough,
there exists unique values
$\bfa(s)$ and $\bfp(s)$ such that Equations
\eqref{eq-monodromy1} and \eqref{eq-monodromy2},
with $\Lambda_i=\Lambda_i(s)$ and $\Pi=\Pi_{\bfa(s),\bfp(s)}$,
are satisfied for
all $(i,j)\in I$.
Moreover, $\bfa(s)$ and $\bfp(s)$ are smooth functions of $s$ for $s\neq 0$, and
extend continuously at $s=0$ with value
$a_{ij}(0)=0$ and $p_{ij}(0)=p_{ij}^0$.
\end{proposition}
\subsection{Choice of an isometry}
\label{section-choice-isometry}
From now on, we fix a couple $(i,j)\in I$. We have in mind to solve Equations
\eqref{eq-monodromy1} and \eqref{eq-monodromy2}.
The computations will be simplified by applying a well chosen isometry.
Let $h$ be an orientation preserving isometry of $\H^3$ such that
$h(S_i)$ is the horosphere $x_3=1$ and $h(S_i\cap S_j)=(0,0,1)$,
in the half-space model.
(This isometry is unique up to composition by a rotation around the vertical axis.)
Since the horospheres $S_i$ and $S_j$ are tangent, $h(S_j)$ is the sphere
of radius $\half$ centered at $(0,0,\half)$.
The limit points of $h(S_i)$ and $h(S_j)$ are respectively $\infty$ and
$0$.

\medskip
As explained in Section \ref{section-isometries},
the isometry $h$ corresponds to a matrix $H\in SL(2,\C)$, unique up to sign.
We have $H\cdot c_i=\infty$ and $H\cdot c_j=0$, where the dot means
the action by homography on the Riemann sphere.
So $H$ may be written in the form
$$H=\frac{1}{\sqrt{c_j-c_i}}\left(\begin{array}{cc}
\rho &-\rho c_j\\
\rho^{-1} & -\rho^{-1} c_i\end{array}\right)$$
where $\rho$ is some complex number.

\medskip

We use hats to denote the action of the isometry $h$ on various objects:
$\widehat{F}=HF$ solves $d\widehat{F}=\whA \widehat{F}$ where
$\whA=H A H^{-1}.$
An elementary computation gives
\begin{equation}
\label{eq-whA}
\whA=\frac{1}{c_j-c_i}\left(\begin{array}{cc}
(G-c_i)(G-c_j) & -\rho^2(G-c_j)^2\\
\rho^{-2}(G-c_i)^2 & -(G-c_i)(G-c_j)\end{array}\right)\myomega.
\end{equation}
The principal solution of $Y'=\whA Y$ is $\whPi= H \Pi H^{-1}$.
Equations \eqref{eq-monodromy1} and \eqref{eq-monodromy2} are equivalent to
\begin{equation}
\label{eq-monodromy1hat}
\whLambda_i^{-1}\whPi(\gamma_{ij})\whLambda_i\in SU(2)
\end{equation}
\begin{equation}
\label{eq-monodromy2hat}
\whLambda_j^{-1}\whPi(\Gamma_{ji})\whLambda_i\in SU(2)
\end{equation}
where
$\whLambda_i=H\Lambda_i$ and
$\whLambda_j=H\Lambda_j$.
\begin{remark} All these quantities: $H$, $\rho$, $\whA$, $\whPi$, $\whLambda_i$,
$\whLambda_j$
actually depend on both indices $i$ and $j$ because
the chosen isometry $h$ does. However, since $i$ and $j$ are fixed until the
very end of Section \ref{section-solution-monodromy}, this dependence will not be written to make notations
lighter.
\end{remark}
\subsection{Computation of the matrices $\whLambda_i(s)$ and $\whLambda_j(s)$}
Consider the matrix
$$\Xi(s)=\left(\begin{array}{cc} e^{s/2} & 0 \\ 0 &e^{-s/2}\end{array}\right)\in SL(2,\C).$$
Then $s\mapsto \Phi(\Xi(s)\Xi(s)^*)$ is the
parametrization at unit speed of the positive vertical axis (oriented upwards)
in the half-space model. (Here $\Phi$ is the isometry from the Minkowski model
to the half-space model given in Section \ref{section-half-space}).
The horosphere $\whS_i=h(S_i)$ is parametrized by $\whf_i=h\circ f_i$.
We need to compute the corresponding null holomorphic map
$\whF_i$.
By substitution of  $G=c_i$ and $\myomega=\lambda_i dz$
in \eqref{eq-whA}, we obtain
\begin{equation}
\label{eq-whAi}
\whA_i=\left(\begin{array}{cc}
0 & \whlambda_i \\ 0 & 0\end{array}\right)
\qquad\qquad
\whlambda_i=\rho^{2}\lambda_i (c_i-c_j)
\end{equation}
By our choice of the isometry $h$, we have $\whF_i(p_{ij}^0)=I_2$.
Hence
\begin{equation}
\label{-eq-whF_i}
\whF_i(z)=\exp((z-p_{ij}^0)\whA_i).
\end{equation}
Then by our choice of the matrix $\Lambda_i(s)$ in Section \ref{section-choice-Lambda},
\begin{equation}
\label{eq-whLambdai}
\whLambda_i(s)=\exp(-p_{ij}^0\whA_i)\Xi(\xi_i s)
\end{equation}
In the same way, the horosphere $\whS_j$ has null holomorphic map $\whF_j$
given by
\begin{equation}
\label{eq-whF_j}
\whF_j(z)=\exp((z-p_{ji}^0)\whA_j).
\end{equation}
\begin{equation}
\label{eq-whAj}
\whA_j=\left(\begin{array}{cc}
0 & 0 \\ \whlambda_j & 0\end{array}\right)
\qquad\qquad
\whlambda_j=\rho^{-2}\lambda_j (c_j-c_i)
\end{equation}
and recalling that the mean curvature vector of $\whS_j$ at $(0,0,1)$ points down,
\begin{equation}
\label{eq-whLambdaj}
\whLambda_j(s)=\exp(-p_{ji}^0\whA_j)\Xi(-\xi_j s)
\end{equation}

\subsection{Partial derivatives of the matrix $\whA$}
We return to the matrix $\whA_{\bfa,\bfp}$ given by \eqref{eq-whA} with
$G=G_{\bft(\bfa)}$ and $\myomega=\myomega_{\bft(\bfa),\bfa}$.
Taking the derivative of \eqref{eq-whA} and using Equations \eqref{eq-G0},
\eqref{eq-omega0}, \eqref{eq-derivee-G}
and \eqref{eq-derivee-omega}, we obtain the following result after simplification:
\begin{proposition}
\label{proposition-whA} At $\bfa=0$, we have:
\begin{equation}
\label{eq-derivee-whAi}
\whA_{0,\bfp}=\whA_i \,dz\qquad\qquad
\frac{\partial\whA_{\bfa,\bfp}}{\partial a_{ij}}=\frac{c_i-c_j}{2}
\left(\begin{array}{cc} 1 & 0 \\ 0 & -1\end{array}\right)\frac{dz}{z-p_{ij}}
\quad\mbox{ in $\C_i$}
\end{equation}
\begin{equation}
\label{eq-derivee-whAj}
\whA_{0,\bfp}=\whA_j \,dz\qquad\qquad
\frac{\partial\whA_{\bfa,\bfp}}{\partial a_{ij}}=\frac{c_j-c_i}{2}
\left(\begin{array}{cc} 1 & 0 \\ 0 & -1\end{array}\right)\frac{dz}{z-p_{ji}}
\quad\mbox{ in $\C_j$}
\end{equation}
\begin{equation}
\label{eq-derivee-whAij}
\whA_{0,\bfp}=0\qquad\qquad
\frac{\partial \whA_{\bfa,\bfp}}{\partial a_{ij}}=
\frac{c_j-c_i}{2}\left(\begin{array}{cc}
z & -\rho^2\\ \rho^{-2}z^2 & -z\end{array}\right)\frac{dz}{z^2}
\qquad\mbox{ in $\CC_{ij}$}
\end{equation}
where the matrices $\whA_i$ and $\whA_j$ are given by \eqref{eq-whAi}
and \eqref{eq-whAj}.
\end{proposition}

\subsection{Expansion of  $\whPi(\gamma_{ij})$}
\begin{proposition}
\label{proposition-whPi-gamma}
We have the following expansion
\begin{equation}
\label{eq-whPi-gamma}
%\whPi_{\bfa,\bfp}(\gamma_{ij})=
%I_2+ a_{ij}\,\pi i (c_i-c_j)\left(\begin{array}{cc} 1 & 2\whlambda_i p_{ij}\\ 0 & -1\end{array}
%\right) + O(||\bfa||^2).
\whPi_{\bfa,\bfp}(\gamma_{ij})=
I_2+ a_{ij}\,\pi \i (c_i-c_j)\exp(-p_{ij}\whA_i)\left(\begin{array}{cc} 1 & 0 \\ 0 & -1\end{array}\right)\exp(p_{ij}\whA_i) + O(||\bfa||^2).
\end{equation}
\end{proposition}
Proof: when $\bfa=0$, the principal solution in $\C_i$ is given by
$\whPi_{0,\bfp}(z,0_i)=\exp(z \whA_i)$, which is well defined, so $\whPi_{0,\bfp}(\gamma_{ij})=I_2$.
By Proposition \ref{proposition-derivative-monodromy} in Appendix \ref{appendix-derivative-monodromy}, Equation \eqref{eq-derivee-whAi} and the residue theorem:
%$$\frac{\partial\whPi_{\bfa,\bfp}(\gamma_{ij})}{\partial a_{ij}}=\frac{c_i-c_j}{2}\int_{\gamma_{ij}}
%\left(\begin{array}{cc}1 & -\whlambda_i z\\ 0 &1\end{array}\right)
%\left(\begin{array}{cc} 1 & 0 \\ 0 & -1\end{array}\right)
%\left(\begin{array}{cc}1 & \whlambda_i z\\ 0 &1\end{array}\right)
%\frac{dz}{z-p_{ij}}.$$
\begin{eqnarray*}
\frac{\partial\whPi_{\bfa,\bfp}(\gamma_{ij})}{\partial a_{ij}}
&=&\frac{c_i-c_j}{2}\int_{\gamma_{ij}}
\exp(-z\whA_i)
\left(\begin{array}{cc} 1 & 0 \\ 0 & -1\end{array}\right)
\exp(z\whA_i)
\frac{dz}{z-p_{ij}}\\
&=&\pi \i (c_i-c_j)\exp(-p_{ij}\whA_i)\left(\begin{array}{cc} 1 & 0 \\ 0 & -1\end{array}\right)\exp(p_{ij}\whA_i).
\end{eqnarray*}
%\begin{equation}
%\frac{\partial\whPi_{\bfa,\bfp}(\gamma_{ij})}{\partial a_{ij}}=
%\pi i (c_i-c_j)\left(\begin{array}{cc} 1 & 2\whlambda_i p_{ij}\\ 0 & -1\end{array}
%\right).
%\end{equation}
For $(k,l)\neq (i,j)$,
the partial derivative of $\whA$ with respect to $a_{k\ell}$ is holomorphic at $p_{ij}$.
Hence the partial derivative of $\whPi(\gamma_{ij})$ with respect to
$a_{k\ell}$ is zero.
Equation \eqref{eq-whPi-gamma} follows.
\cqfd

\subsection{Expansion of $\whPi(\Gamma_{ji})$.}
Recall that $\exp$ is a local diffeomorphism from a neighborhood of 
$0$ in $\boM_2(\C)$ to a neighborhood of $I_2$ in $GL(2,\C)$.
We denote the inverse diffeomorphism by $\log$. Also, for a matrix $M$ and
a complex number $\lambda$, we define $M^{\lambda}=\exp(\lambda\log M)$.
\begin{proposition}
\label{proposition-whPi-Gamma}
We have
\begin{equation}
\label{eq-whPi-Gamma}
\whPi_{\bfa,\bfp}(\Gamma_{ji})\times\whPi_{\bfa,\bfp}(\gamma_{ij})^{\frac{-2}{2\pi \i}\log a_{ij}}=
\exp(-p_{ji}\whA_j)\exp(p_{ij}\whA_i)+{\mathcal O}(\bfa)
\end{equation}
where ${\mathcal O}(\bfa)$ is a well defined holomorphic function of $(\bfa,\bfp)$ which vanishes at $\bfa=0$.
\end{proposition}
Proof: Let us fix the value of all parameters except $a_{ij}$. We must first see that the left-hand side of \eqref{eq-whPi-Gamma} is a well defined
function of $a_{ij}$ for $a_{ij}\neq 0$, which of course means that $\whPi(\Gamma_{ji})$
itself is not. To see this, we have to define precisely the path $\Gamma_{ji}$ from
$0_i$ to $0_j$.
For $t_{ij}$ and $t_{ji}$ non zero, we define $\Gamma_{ji}$  as the composition of the following
two paths:
\begin{itemize}
\item The following path from $0_i$ to $p_{ij}+t_{ij}$ in $\C_i$:
a path from $0_i$ to $p_{ij}+1$ in $\C_i$ minus all unit disks around the nodes
(depending continuously on the parameter $p_{ij}$ in a neighborhood of $p_{ij}^0$),
composed with the spiral parametrized by $x\mapsto p_{ij}+(t_{ij})^x$ for $x\in [0,1]$.
\item The following path from $p_{ji}+t_{ji}$ to $0_j$ in $\C_j$:
the spiral parametrized by $x\mapsto p_{ji}+(t_{ji})^{1-x}$ for $x\in [0,1]$, composed with a path from $p_{ji}+1$ to $0_j$.
\end{itemize}
We can compose these two paths because the points $p_{ij}+t_{ij}$ and
$p_{ji}+t_{ji}$ are both identified with $1_{ij}$ when opening nodes.
Observe that we need a determination of the arguments of $t_{ij}$ and $t_{ji}$
to define the spirals. In other words, if we take $t_{ij}$ and $t_{ji}$
to live in the universal cover of the punctured unit disk (so $\arg t_{ij}$
and $\arg t_{ji}$ are well defined), $\Gamma_{ji}$ depends continuously on  $t_{ij}$ and $t_{ji}$.
Now if the argument
of $a_{ij}$ is increased by $2\pi$, then the arguments of
$t_{ij}$ and $t_{ji}$ are increased by the same amount. Hence the homotopy class
of $\Gamma_{ji}$
is multiplied on the right by $(\gamma_{ij})^2$
and $\whPi(\Gamma_{ji})$ is multiplied on the right by $\whPi(\gamma_{ij})^2$.
Consequently, the left hand side of \eqref{eq-whPi-Gamma} is unchanged, so
is a well defined holomorphic function of $a_{ij}$ for $a_{ij}\neq 0$.

\medskip

Next we prove that the left hand side of \eqref{eq-whPi-Gamma} is uniformly bounded. Because it is a well defined
function of $a_{ij}$, we can assume that $\arg a_{ij}\in[-2\pi,2\pi]$.
Using \eqref{eq-morphism}, we write
\begin{equation}
\label{eq-Pi-produit}
\whPi(\Gamma_{ji})=\whPi(0_j,p_{ji}+1)\whPi(p_{ji}+1,p_{ji}+t_{ji})
\whPi(p_{ij}+t_{ij},p_{ij}+1)\whPi(p_{ij}+1,0_i).
\end{equation}
Since the path from $0_i$ to $p_{ij}+1$ stays in a fixed compact set of
$\C_i$ minus the nodes, where $\whA$ is uniformly bounded, the
fourth factor in \eqref{eq-Pi-produit} is uniformly bounded.
We estimate the third factor using Proposition \ref{proposition-annulus} from
Appendix \ref{appendix-annulus}.
For this, we need an integral estimate of $||\whA||$ on the circle
of center $p_{ij}$ and radius $|t_{ij}|/2$.

\medskip

I claim that $\whA=O(a_{ij})$ in compact subsets of $\CC_{ij}\setminus\{0,\infty\}$
(even if the other parameters $a_{k\ell}$ are non-zero).
Indeed, $\whA$ depends holomorphically on $a_{ij}$, and if $a_{ij}=0$, then
by \eqref{eq-aij-nul}, $t_{ij}=t_{ji}=0$, so $\Omega=0$ and $\whA=0$ in $\CC_{ij}$.
Also, by \eqref{eq-derivee-tij}, $a_{ij}=O(t_{ij})$.
Consequently, since $\whA$ is a matrix-valued 1-form,
$$
\int_{|z-p_{ij}|=\frac{|t_{ij}|}{2}}||\whA||=\int_{|v_{ij}|=\frac{|t_{ij}|}{2}}||\whA||
=\int_{|w_{ij}|=2}||\whA||\leq C|t_{ij}|$$
for some uniform constant $C$.
By Proposition \ref{proposition-annulus} in Appendix \ref{appendix-annulus},
$\whPi(p_{ij}+t_{ij},p_{ij}+1)$ is uniformly bounded.
The first and second factors in \eqref{eq-Pi-produit} are estimated in the 
exact same way.
We conclude that
$\whPi(\Gamma_{ji})$ is uniformly bounded (although not well defined --
but we assumed that $\arg a_{ij}\in[-2\pi,2\pi]$).
The left-hand side of
\eqref{eq-whPi-Gamma} is now a bounded, well defined holomorphic function of 
$(\bfa,\bfp)$ on the set $a_{ij}\neq 0$.
By Riemann extension theorem (in several variables),
it extends holomorphically at $a_{ij}=0$.

\medskip

To compute its value at $\bfa=0$, assume that all parameters $a_{k\ell}$ for
$(k,\ell)\neq (i,j)$ are zero.
By \eqref{eq-derivee-whAi}, $\whA_{\bfa,\bfp}-\whA_i=O(a_{ij})$ in compact subsets of $\C_i$ minus the nodes.
By Point (2) of Proposition \ref{proposition-annulus} (with $\wtA=\whA_i$), we obtain
$$||\whPi_{\bfa,\bfp}(p_{ij}+t_{ij},0_i)-\whPi_i(p_{ij}+t_{ij},0_i)||\leq C\left|a_{ij}\log |a_{ij}|\right|$$
where $\whPi_i$ is the principal solution of $Y'=\whA_i Y$ in $\CC_i$, namely
$\whPi_i(z,0_i)=\exp(z \whA_i).$
This gives
$$\lim_{a_{ij}\to 0} \whPi_{\bfa,\bfp}(p_{ij}+t_{ij},0_i)=\exp(p_{ij}\whA_i ).$$
Arguing in the same way, we obtain
$$\lim_{a_{ij}\to 0}\whPi_{\bfa,\bfp}(0_j,p_{ji}+t_{ji})=\exp(-p_{ji}\whA_j ).$$
Proposition \ref{proposition-whPi-Gamma} follows.
\cqfd
\subsection{Solution of the monodromy problem}
\label{section-solution-monodromy}
We are now ready to prove Proposition \ref{proposition-monodromy}.
The unitary group $SU(2)$ is not a complex manifold so we have to leave the
realm of holomorphic functions.
We introduce a small positive real number $\tau$ and have in mind to apply
the implicit function theorem at $\tau=0$. We write
$$a_{ij}=\tau\frac{b_{ij}}{c_i-c_j}$$
where $b_{ij}$ is a complex number in a neighborhood of a non-zero central value
$b_{ij}^0$.
The computation will be simplified by knowing a priori the order of
each parameter as a function of $\tau$. The correct orders are
$$s=-\tau\log\tau$$
$$p_{ij}=p_{ij}^0+sq_{ij}$$
where $q_{ij}$ is a complex parameter in a neighborhood of $0$.
One issue here is that the function $\tau\mapsto \tau\log\tau$ does not extend
as a differentiable function at $\tau=0$. We solve this problem by writing
$\tau=e^{-1/t^2}$
where $t$ is a real parameter in a neighborhood of $0$.
Both $\tau$ and $\tau\log\tau$ extend smoothly at $t=0$, and all parameters
are smooth functions of $t$.
Let $\bfb=(b_{ij})_{(i,j)\in I}$ and $\bfq=(q_{ij},q_{ji})_{(i,j)\in I}$.

\medskip
Recall that $\exp$ maps
the Lie algebras $\mathfrak{sl}(2,\C)$ and $\mathfrak{su}(2,\C)$ to
the Lie groups $SL(2,\C)$ and $SU(2,\C)$, respectively.
We define
$$P_{ij}=P_{ij}(t,\bfb,\bfq)=\log\left(\whLambda_i(s)^{-1}\whPi_{\bfa,\bfp}(\gamma_{ij})\whLambda_i(s)\right)
\in\mathfrak{sl}(2,\C)$$
$$Q_{ij}=Q_{ij}(t,\bfb,\bfq)=\log\left(\whLambda_j(s)^{-1}\whPi_{\bfa,\bfp}(\Gamma_{ji})\whLambda_i(s)\right)
\in\mathfrak{sl}(2,\C).$$
We want to solve $P_{ij}\in \mathfrak{su}(2,\C)$ and $Q_{ij}\in \mathfrak{su}(2,\C)$.
We compute $P_{ij}$ using Proposition \ref{proposition-whPi-gamma}:
\begin{equation}
\label{eq-solution1}
\whPi(\gamma_{ij})=I_2+\pi \i \tau b_{ij}\exp(-p_{ij}\whA_i)
\left(\begin{array}{cc} 1 & 0 \\
0 & -1\end{array}\right)\exp(p_{ij}\whA_i)+O(\tau^2)
\qquad\mbox{ using \eqref{eq-whPi-gamma}}.
\end{equation}

$$P_{ij}=
\pi \i\tau b_{ij}\,\Xi(-\xi_i s)\exp(-sq_{ij}\whA_i)\left(\begin{array}{cc} 1 & 0 \\
0 & -1\end{array}\right)\exp(sq_{ij}\whA_i)\Xi(\xi_i s)+O(\tau^2)
\qquad\mbox{ using \eqref{eq-whLambdai}}$$
\begin{equation}
\label{eq-Pij}
P_{ij}=\pi \i\tau b_{ij}
\left(\begin{array}{cc} 1 & 2\whlambda_is q_{ij}\\ 0 & -1\end{array}\right)
+O(\tau^2)
\end{equation}
We compute $Q_{ij}$ using Proposition \ref{proposition-whPi-Gamma}:
$$\whPi(\Gamma_{ji})=\exp(-p_{ji}\whA_j)\exp(p_{ij}\whA_i)
\whPi(\gamma_{ij})^{\frac{2}{2\pi \i}\log a_{ij}}+O(\tau)
\qquad\mbox{ using \eqref{eq-whPi-Gamma}}.$$
$$\whPi(\gamma_{ij})^{\frac{2}{2\pi \i}\log a_{ij}}=
I_2-s b_{ij}\exp(-p_{ij}\whA_i)\left(\begin{array}{cc} 1 & 0 \\
0 & -1\end{array}\right)\exp(p_{ij}\whA_i)+O(\tau)
\qquad\mbox{ using \eqref{eq-solution1}}.$$
$$
\whLambda_j(s)^{-1}\whPi(\Gamma_{ji})\whLambda_i(s)
=
\Xi(\xi_j s)\exp(-sq_{ji}\whA_j)\left(\begin{array}{cc}
1-s b_{ij} & 0 \\ 0 &1+s b_{ij}\end{array}\right)\exp(sq_{ij}\whA_i)\Xi(\xi_i s)+O(\tau).$$
\begin{equation}
\label{eq-Qij}
Q_{ij}=s\left(\begin{array}{cc} \frac{\xi_i+\xi_j}{2}-b_{ij}& \whlambda_i q_{ij}\\ -\whlambda_j q_{ji}&  b_{ij}-\frac{\xi_i+\xi_j}{2}\end{array}\right)+O(\tau).
\end{equation}
For a matrix $M=(M^{k\ell})_{1\leq k,\ell\leq 2}\in\mathfrak{sl}(2,\C)$, we have
$$M\in\mathfrak{su}(2,\C)\quad\Leftrightarrow\quad
\Re(M^{11})=0\quad\mbox{ and } \quad M^{12}+\overline{M^{21}}=0.$$
Define the function $\boF=(\boF_{ij})_{(i,j)\in I}$ for $t\neq 0$ by
$$\boF_{ij}(t,\bfb,\bfq)=
\left(\frac{1}{\tau}\Re(P_{ij}^{11})\;,\;
\frac{1}{\tau s}(P_{ij}^{12}+\overline{P_{ij}^{21}})\;,\;
\frac{1}{s}\Re(Q_{ij}^{11})\;,\;
\frac{1}{s}(Q_{ij}^{12}+\overline{Q_{ij}^{21}})\right).$$
We want to solve $\boF(t,\bfb,\bfq)=0$.
By \eqref{eq-Pij} and \eqref{eq-Qij}, $\boF$ extends smoothly at $t=0$, with
$$\boF_{ij}(0,\bfb,\bfq)=
\left(-\pi\Im(b_{ij}) \;,\;
2\pi \i \,b_{ij}\whlambda_i q_{ij}
\;,\;\mbox{$\frac{\xi_i+\xi_j}{2}$}-\Re(b_{ij})
\;,\;
\whlambda_iq_{ij}-\overline{\whlambda_jq_{ji}}\;\right).$$
Taking $b^0_{ij}=\frac{\xi_i+\xi_j}{2}>0$, we have
$\boF(0,\bfb^0,0)=0.$
It is straightforward that the partial differential of $\boF$ with respect to
the variables $(\bfb,\bfq)$ at $(0,\bfb^0,0)$ is an  isomorphism.
By the implicit function theorem, for $t$ in a neighborhood of $0$,
there exists $(\bfb(t),\bfq(t))$ depending smoothly on $t$ such that
$\boF(t,\bfb(t),\bfq(t))=0$. Proposition \ref{proposition-monodromy} is proved, and the monodromy
problem is solved.
\cqfd
\section{Embeddedness}
\label{section-embedded}
Here is what we have achieved so far. For each small enough value of the parameter $t>0$,
we have constructed a null holomorphic map $F$ which has $SU(2)$-valued monodromy.
All parameters are now smooth functions of $t>0$.
To ease notation, the dependence on $t$ will not be written.
Let $f:\Sigma\to\H^3$ be the CMC-1 immersion associated to $F$.
It remains to prove that $f(\Sigma)$ is embedded.
We work in the half-space model, so $f(z)=\Phi(F(z)F(z)^*)$ is given by formula \eqref{eq-half-space}.
Fix a small number $\varepsilon>0$. We consider the following disjoint domains in
$\Sigma$:
$$\C_i^{\varepsilon}=\{z\in\C_i\,:\, \forall j\in J_i,\;|z-p_{ij}^0|>\varepsilon\}$$
$$\CC_{ij}^{\varepsilon}=\{z\in \CC_{ij}\,:\,\varepsilon<|z|<\frac{1}{\varepsilon}\}.$$
The complement of these domains in $\Sigma$ are annuli which we call transition regions.
We also fix some large number $R$ and define $\C_i^{\varepsilon,R}=\C_i\cap
D(0,R)$.
\subsection{Geometry of the image of $\C_i^{\varepsilon,R}$}
\label{section-embedded1}
Fix some $i$, $1\leq i\leq n$, and consider an isometry $h$ such that
$h(S_i)$ is the horosphere $x_3=1$ and $h$ maps $f_i(0)$ to the point $(0,0,1)$. The isometry $h$ is represented by a matrix $H\in SL(2,\C)$
which has the form
$$H=\frac{1}{\sqrt{c-c_i}}\left(\begin{array}{cc} \rho & -\rho c\\\rho^{-1} & -\rho^{-1} c_i\end{array}
\right)$$
where $\rho, c$ are some complex numbers ($\rho$ not the same as in Section
\ref{section-choice-isometry}).
As in Section \ref{section-choice-isometry}, we use hats to denote the action of $h$,
so $\whf=h\circ f$, $\whF=HF$ and so on.
Equations \eqref{eq-whA} and \eqref{eq-whAi} hold true, with $c$ in place
of $c_j$.
By construction, $\whf(0_i)$ parametrizes the vertical axis at speed
$\xi_i$ as $s$ varies, so $\whF(0_i)=\Xi(\xi_i s)$, up to right multiplication by
$SU(2)$.
We have $\whA=\whA_i+O(\tau)$. Since $\C_i^{\varepsilon,R}$ is a fixed compact domain,
$$\whF(z)=\exp(z\whA_i)\Xi(\xi_i s)+O(\tau) \qquad \mbox{ for $z\in \C_i^{\varepsilon,R}$}.$$
From this, we conclude that $\whf(\C_i^{\varepsilon,R})$ converges smoothly to (a subdomain of) the
horosphere $x_3=1$ as $t\to 0$.
Moreover, from \eqref{eq-half-space}, we get
$$x_3(z)=e^{\xi_i s}+O(\tau) \qquad \mbox{ for $z\in \C_i^{\varepsilon,R}$}$$
so for $t$ small enough, the image of $\C_i^{\varepsilon,R}$
lies above the horosphere $x_3=1$.
\subsection{Geometry of the end at $\infty_i$}
\label{section-embedded2}
\def\sommexi{{\zeta}}
Next we prove that the image of $|z|>R$ in $\C_i$ is embedded.
I claim that for $t>0$ small enough, the Gauss map $G$ has multiplicity $1$ at $\infty_i$. This is delicate because $G$ is constant when $t=0$.
We work in the local coordinate $w=\frac{1}{z}$ in a neighborhood of $\infty_i$
and write $\wtG(w)=G(1/w)$.
From \eqref{eq-derivee-G}, we obtain
$$\frac{\partial \wtG'(w)}{\partial a_{ij}}=\frac{c_j-c_i}{2\lambda_i}\frac{1}{(1-p_{ij}w)^2}
\qquad\mbox{ for $j\in J_i^+$}$$
$$\frac{\partial \wtG'(w)}{\partial a_{ji}}=\frac{c_i-c_j}{2\lambda_i}\frac{1}{(1-p_{ji}w)^2}
\qquad\mbox{ for $j\in J_i^-$}.$$
\begin{eqnarray*}
\wtG'(0)&=&\sum_{j\in J_i^+}\frac{\partial \wtG'(0)}{\partial a_{ij}} a_{ij}+
\sum_{j\in J_i^-}\frac{\partial \wtG'(0)}{\partial a_{ji}} a_{ji} + O(||\bfa||^2)\\
&=& \sum_{j\in J_i^+} \frac{c_j-c_i}{2\lambda_i}\frac{\tau b_{ij}}{c_i-c_j}
+\sum_{j\in J_i^-} \frac{c_i-c_j}{2\lambda_i}\frac{\tau b_{ji}}{c_j-c_i}+O(\tau^2)\\
&=&-\frac{\tau\sommexi_i}{2\lambda_i}+o(\tau)
\qquad\mbox{ where} \qquad\sommexi_i=\frac{1}{2}\sum_{j\in J_i}(\xi_i+\xi_j)>0.
\end{eqnarray*}
Hence for $t>0$ small enough, $\wtG'(0)\neq 0$, so the Gauss map has multiplicity
one at the end.
To study the geometry of the end, we consider again the isometry $h$ introduced in
Section \ref{section-embedded1}. Then $\whG=H\cdot \wtG$ has a simple pole at $w=0$
with residue

$$\Res_{w=0} \;\whG=\Res_{w=0}\;\rho^2\frac{\wtG-c}{\wtG-c_i}
=\rho^2\frac{c_i-c}{\wtG'(0)}
\simeq\frac{-2\whlambda_i}{\tau\sommexi_i}$$
where
$\whlambda_i=\rho^2\lambda_i(c_i-c).$
From \eqref{eq-whA} we obtain
$$\whG^2\widehat{\Omega}=-\whA_{21}\simeq -\whlambda_i dz=\whlambda_i\frac{dw}{w^2}.$$
By Theorem \ref{theorem-end} in Appendix \ref{appendix-end} (with
$\alpha=\frac{-2\whlambda_i}{\tau\zeta_i}$ and
$\alpha^2\beta=\whlambda_i$), there exists a
\underline{uniform} positive $\epsilon$ (independent of $t$) such that the image of $0<|w|<\epsilon$ is the vertical graph $x_3=u(x_1,x_2)$ of a function $u$. Moreover, at infinity
we have
$$\log u(x_1,x_2)\simeq (\tau\sommexi_i+o(\tau))\log\sqrt{x_1^2+x_2^2}$$
 so $x_3>1$ on the end.
Replacing $R$ by $\epsilon^{-1}$ if necessary, we obtain that
$\whf(\C_i^{\varepsilon})$ is embedded.
and moreover lies above the horosphere $x_3=1$ (using the maximum
principle).
In other words, $f(\C_i^{\varepsilon})$ lies on the mean-convex side of the horosphere $S_i$.
This guarantees that the images $f(\C_i^{\varepsilon})$ for $1\leq i\leq n$ are disjoint.
\begin{remark}
\label{remark-ignore}
From this, we conclude that we can always ignore a tangency point
by simply removing the corresponding couple $(i,j)$ from $I$, and still obtain an
embedded surface for $t>0$.
\end{remark}
\subsection{Geometry of the catenoidal necks.}
Fix a couple $(i,j)\in I$. Consider again the isometry $h$ introduced in 
Section \ref{section-choice-isometry}, which maps the horosphere $S_i$
to the horosphere $x_3=1$ and the horosphere $S_j$ to the
sphere of radius $\frac{1}{2}$ centered at $(0,0,\frac{1}{2})$.
In this section, we prove that after a blowup of ratio $1/\tau$,
the image $\whf(\CC_{ij}^{\varepsilon})$ converges to a vertical catenoid.

\medskip

By computations similar to the computation of $Q_{ij}$ in Section \ref{section-solution-monodromy}, we have
\begin{equation}
\label{eq-catenoid1}
\whF(1_{ij})=\whPi(p_{ij}+t_{ij},O_i)\whF(0_i)=
I_2+\frac{s}{2}\left(\begin{array}{cc} \xi_i-b_{ij} & 2\whlambda_i q_{ij}\\
0 & b_{ij}-\xi_i\end{array}\right)+O(\tau).
\end{equation}
By \eqref{eq-derivee-whAij}, we have in $\CC_{ij}^{\varepsilon}$
\begin{equation}
\label{eq-catenoid2}
\whA(z)=\tau\wtA(z)+O(\tau^2)\quad\mbox{ with }\quad
\wtA(z)=-\frac{b_{ij}}{2}\left(\begin{array}{cc}z & -\rho^2\\ \rho^{-2} z^2 & -z\end{array}\right)\frac{dz}{z^2}.
\end{equation}
Let
$$\wtF(z)=\frac{1}{\tau}(\whF(z)-\whF(1_{ij})).$$
Then using \eqref{eq-catenoid1} and \eqref{eq-catenoid2},
$$d\wtF(z)=\frac{1}{\tau}\whA(z)\whF(z)
=(\wtA(z)+O(\tau))(I_2+O(s))=\wtA(z)+O(s).$$
Since $\CC_{ij}^{\varepsilon}$ is a fixed compact set, we obtain by integration
$$\wtF(z)=-\frac{b_{ij}}{2}\left(\begin{array}{cc}
\log z & \rho^2(z^{-1}-1) \\ \rho^{-2} (z-1) & -\log z\end{array}\right)+O(s).$$
Write $\whx_k(z)=\whx_k(1_{ij})+\tau \wtx_k(z)$ for $1\leq k\leq 3$.
Using \eqref{eq-half-space} and $\whF(z)=I_2+O(s)$, we obtain
$$\left(\wtx_1(z)+\i\wtx_2(z),\wtx_3(z)\right)=
\left(\wtF_{12}(z)+\overline{\wtF_{21}(z)},
-2\,\Re\left(\wtF_{22}(z)\right)\right)+O(s).$$
$$\lim_{t\to 0}\left(\wtx_1(z)+\i\wtx_2(z),\wtx_3(z)\right)=-\frac{\xi_i+\xi_j}{4}\left(
\rho^2(z^{-1}-1)+\overline{\rho^{-2}(z-1)},2\log|z|\right).$$
This is the parametrization of a vertical catenoid of necksize $\frac{\xi_i+\xi_j}{2}$.
This means that after a blowup of ratio $\frac{1}{\tau}$ at $\whf(1_{ij})$, the image of $\CC_{ij}^{\varepsilon}$ converges smoothly to a catenoid. Also observe that the image
of the circle $|z|=\varepsilon$ lies above the image of $|z|=\frac{1}{\varepsilon}$.
Finally, \eqref{eq-catenoid1} gives
$$\whx_3(1_{ij})=1+s\frac{\xi_i-\xi_j}{2}+O(\tau).$$
Hence the catenoidal neck lies below the image of $\C_i^{\varepsilon}$.
\subsection{Geometry of the transition regions}
Fix $(i,j)\in I$ and let $U_{ij}$ be the annulus in $\Sigma$ bounded by the circles
$|z-p_{ij}|=\varepsilon$ in $\C_i$ and $|z|=\varepsilon$ in $\CC_{ij}$.
We consider again the isometry $h$ introduced in Section \ref{section-choice-isometry}.
Let us prove that the mean curvature vector of
$\whf$ is almost vertical in $U_{ij}$. Given the geometric interpretation of the
Gauss map given in Section \ref{section-half-space}, an elementary
computation shows that the angle $\theta(z)$ between the mean curvature
vector at $f(z)$ and the vertical axis is related to the Gauss map $G(z)$ by
\begin{equation}
\label{eq-theta-G}
\frac{\sin\theta(z)}{1+\cos\theta(z)}=\frac{x_3(z)}{|G(z)-x_1(z)-\i x_2(z)|}.
\end{equation}
The function $\whG^{-1}$ is holomorphic in the annulus $U_{ij}$ and is bounded
by $C\varepsilon$ on the boundary circles, for some
uniform constant $C$.
By the maximum principle, $\whG^{-1}$ is bounded by $C\varepsilon$ in $U_{ij}$.
The norm of the holomorphic map $F(z)-I_2$ is bounded by $C\varepsilon$ on the
boundary of $U_{ij}$, so is bounded by $C\varepsilon$ in $U_{ij}$ by the maximum
principle. Hence
the function $\whx_1+\i\whx_2$ is uniformly bounded in $U_{ij}$, and
the height $\whx_3$ satisfies $|\whx_3-1|\leq C\varepsilon$ in $U_{ij}$.
Using \eqref{eq-theta-G}, we obtain
$$\frac{\sin\widehat{\theta}(z)}{1+\cos\widehat{\theta}(z)}\leq C\varepsilon \quad\mbox{in $U_{ij}$}.$$
Hence by choosing $\varepsilon$ small enough, we can ensure that $\widehat{\theta}(z)<\frac{\pi}{2}$.
This implies that $\whf(U_{ij})$
is locally a vertical graph. Since we have already seen that it is a graph on the boundary circles, it is globally a graph.
Moreover, by the maximum principle, it lies above the lowest point of the top boundary component of the catenoidal neck $\whf(\CC_{ij}^{\varepsilon})$.
The image of the annulus bounded by the circles
$|z-p_{ji}|=\varepsilon$ in $\C_j$ and $|z|=\frac{1}{\varepsilon}$ in $\CC_{ij}$
is studied in the same way, using an isometry which maps the horosphere $S_j$
to the horosphere $x_3=1$.
This proves that $f(\Sigma)$ is embedded and concludes the proof of Theorem \ref{main-theorem}.\cqfd

\appendix
\section{Derivative of the monodromy}
\label{appendix-derivative-monodromy}
Consider a domain $\Omega\subset\C$ and a point $z_0\in\Omega$.
Let $A_{\lambda}(z)\in GL(n,\C)$
be a family of matrices depending holomorphically on $(\lambda,z)$ for 
$z\in\Omega$ and $\lambda$ in a neighborhood of $0$.
Let $\Pi_{\lambda}$ denote the principal solution of $Y'=A_{\lambda} Y$
in $\Omega$.
\begin{proposition}
\label{proposition-derivative-monodromy}
For any $\gamma\in\pi_1(\Omega,z_0)$,
$$\frac{\partial\Pi_{\lambda}(\gamma)}{\partial\lambda}|_{\lambda=0}
=\Pi_0(\gamma)\int_{\gamma}\Pi_0(z,z_0)^{-1}\frac{\partial A_{\lambda}(z)}{\partial\lambda}
\Pi_0(z,z_0)\,dz.$$
\end{proposition}
Proof: Let $Y_{\lambda}(z)=\Pi_{\lambda}(z,z_0)$ and
$W=\partial Y_{\lambda}/\partial\lambda$.
Differentiating $Y_{\lambda}'=A_{\lambda}Y_{\lambda}$ and $Y_{\lambda}(z_0)=I_n$
 with respect to $\lambda$
at $\lambda=0$, we get $W(z_0)=0$ and
$$W'=A_0 W+\frac{\partial A_{\lambda}}{\partial\lambda} Y_0.$$
By the variation of constants formula (Theorem 3.12 in \cite{teschl})
$$W(z)=Y_0(z)\int_{z_0}^z Y_0(w)^{-1}\frac{\partial A_{\lambda}(w)}{\partial\lambda}
Y_0(w)\,dw.$$
Taking $z=\gamma(1)$, the result follows.
\cqfd
\section{Uniform estimates of the solution of $Y'=AY$ in an annulus}
\label{appendix-annulus}
In this section, we consider the annulus $\Omega\subset\C$ defined by
$\rho^{-1}t < |z| < \rho$, where $\rho>1$ is some
fixed number and $t$ is a small positive parameter. We are aiming for estimates
which are uniform with respect to $t$.
Let $A:\Omega\to \mathfrak{sl}(n,\C)$ be a holomorphic map.
Let $Y(z)\in SL(n,\C)$ be the solution of $Y'=A Y$ in $\Omega$, with initial condition $Y(1)=I_n$.
(Of course, $Y(z)$ is only well defined in the universal cover of $\Omega$:
the value of $Y(z)$ depends on the determination of $\arg z$.)

\begin{proposition}
\label{proposition-annulus}
\begin{enumerate}
\item Assume that for some constant $c$,
$$
\int_{|z|=\rho}||A||\leq c
\quad\mbox{ and } \quad
\int_{|z|=\rho^{-1}t} ||A||\leq c t.
$$
Then for $t\leq |z|\leq 1$ and $|\arg z|\leq c'$, $||Y(z)||$ is bounded by a constant
depending only on $c$, $c'$ and $\rho$.
\item Let $\wtA$ be another matrix-valued map satisfying the same hypotheses as $A$ and let $\wtY(z)$ be the solution of $\wtY'=\wtA \wtY$ with initial condition $\wtY(1)=I_n$.
Assume moreover that
$$
\int_{|z|=\rho}||A-\wtA||\leq c t.
$$
Then for $t\leq |z|\leq 1$ and $|\arg(z)|\leq c'$,
$$||Y(z)-\wtY(z)||\leq Ct|\log t|$$
for some constant $C$ depending only on $c$, $c'$ and $\rho$.
\end{enumerate}
\end{proposition}
Proof. We use the letter $C$ for uniform constants, depending only
on $c$ and $\rho$ but not on $t$.
By Cauchy theorem,
$$A(z)=\frac{1}{2\pi \i}\int_{|w|=\rho}\frac{A(w)}{w-z}dw-\frac{1}{2\pi \i}\int_{|w|=\rho^{-1}t}\frac{A(w)}{w-z}dw.$$
Hence for $t\leq |z|\leq 1$,
$$||A(z)||\leq \frac{1}{2\pi}\left(\frac{c}{\rho-1}+ \frac{ct}{t(1-\rho^{-1})}\right)
\leq C.$$
We can connect $1$ and $z$ (in the universal cover of $\Omega$) by a path of length
less than $1+c'$. 
The first Point of Proposition \ref{proposition-annulus} follows from Gromwall
inequality (Lemma 2.7 in \cite{teschl}).
Using Cauchy formula in the same way, we obtain
for $t\leq |z|\leq 1$
$$||A(z)-\wtA(z)||\leq \frac{1}{2\pi}\left(\frac{ct}{\rho-1}+ \frac{ct}{|z|-\rho^{-1}t}\right)
\leq \frac{Ct}{|z|}.$$
By the variation of constants formula (Theorem 3.12 in \cite{teschl})
$$\wtY(z)=Y(z)+Y(z)\int_1^z Y(w)^{-1}(\wtA(w)-A(w))\wtY(w)dw.$$
This gives
$$||\wtY(z)-Y(z)||\leq C t\int_1^z\frac{|dw|}{|w|}
\leq Ct\left(|\arg z|+|\log t|\right).$$
\cqfd
\section{Embedded CMC-1 ends}
\label{appendix-end}
\begin{theorem}
\label{theorem-end}
Let $f:D^*(0,1)\to \H^3$ be a conformal, CMC-1 immersion of the punctured closed unit disk.
Assume that the Gauss map $G$ has a simple pole at $0$, with residue $\alpha$, and
the holomorphic differential $\Omega$ is holomorphic at $0$, with
$\Omega(0)=\beta dz$. Assume that $0<|\alpha\beta|\leq \frac{1}{8}$.
Then there exists $\varepsilon>0$ such that in the half-space model,
$f(D^*(0,\varepsilon))$ is a vertical graph
$x_3=u(x_1,x_2)$ of a (positive) function $u$ over an
exterior domain in the plane. 
The number $\varepsilon$ only depends on a bound on $|\alpha^2\beta|^{\pm 1}$
and $||F(z)||$ on the unit circle.
Moreover, $\alpha\beta$ is real and at infinity, the function $u$ has the following asymptotic behavior:
$$\log u(x_1,x_2)\simeq (1-\sqrt{1+4\alpha\beta})\log |x_1+\i x_2|.$$
\end{theorem}
\begin{remark}
In this paper, we are interested in the case where $\alpha\beta\to 0$ and we
need a uniform positive $\varepsilon$. The conclusions of Theorem
\ref{theorem-end} remain true without the hypothesis $|\alpha\beta|\leq \frac{1}{8}$
but the proof is more involved, as the fuchsian system can be resonant.
In particular, one can prove that $\alpha\beta$ is always a real number
in $(-\frac{1}{4},\infty)$.
\end{remark}
Proof: We use the theory of fuchsian systems to compute $F(z)$ such that
$f=FF^*$ in the punctured disk.
The system $F'=AF$ is fuchsian provided the matrix $A(z)$ has a simple pole
at $0$, which is not the case here (it has a double pole). To circumvent this problem, we introduce the matrix
\begin{equation}
\label{eq-N}
N(z)=\left(\begin{array}{cc} 1 & 0 \\ 0 & z\end{array}\right)
\end{equation}
and make the change of unknown
$F(z)=N(z) \wtF(z).$
By a straightforward computation, \eqref{eq-F} is equivalent to
\begin{equation}
\label{eq-wtF}
\wtF'(z)=\wtA(z)\wtF(z)
\end{equation}
where
$$\wtA=\left(\begin{array}{cc}
G \omega & -z G^2\omega\\  z^{-1}\omega & -G\omega-z^{-1}
\end{array}\right).$$
Now the matrix $\wtA$ has a simple pole at $0$, with residue
$$A_0=\Res_0\wtA=\left(\begin{array}{cc}
\alpha\beta & -\alpha^2\beta \\ \beta & -\alpha\beta -1
\end{array}\right)$$
and the system \eqref{eq-wtF} is fuchsian.
The eigenvalues of $A_0$ are
$$\lambda_1=\frac{-1+\sqrt{\Delta}}{2},\quad
\lambda_2=\frac{-1-\sqrt{\Delta}}{2}
\qquad \mbox{ where }
\Delta=1+4\alpha\beta.$$
The system \eqref{eq-wtF} is called resonant if $\lambda_1-\lambda_2=\sqrt{\Delta}$ is a non-zero integer. It follows from our hypothesis that $\frac{1}{2}\leq |\Delta|\leq \frac{3}{2}$ and $\Delta\neq 1$, so the system is non-resonant.
By the standard theory of fuchsian systems (Proposition 11.2 in \cite{taylor}),
the solution of \eqref{eq-wtF} has the form
\begin{equation}
\label{eq-solution-wtF}
\wtF(z)=U(z) z^{A_0} Y_0
\end{equation}
where $U(z)\in GL(2,\C)$ is well defined, holomorphic in $D(0,1)$ and satisfies $U(0)=I_2$, and $Y_0\in GL(2,\C)$ is a constant matrix.
The monodromy of $\wtF$ on the unit circle $\gamma$ is
$$M_{\gamma}(\wtF)=M_{\gamma}(F)=Y_0^{-1} e^{2\pi \i A_0} Y_0.$$
Since $f$ is well defined, its monodromy $M_{\gamma}(F)$ belongs to $SU(2)$ so its eigenvalues are complex number of modulus 1. This implies
that the eigenvalues $\lambda_1$, $\lambda_2$ of $A_0$ are real numbers,
so $\alpha\beta$ is real.
To compute $z^{A_0}$, we write $A_0=PDP^{-1}$ with
$$P=\frac{1}{\Delta^{1/4}}\left(\begin{array}{cc}
-\lambda_2 & \alpha\lambda_1\\\alpha^{-1} \lambda_1 & -\lambda_2
\end{array}\right)\in SL(2,\C),\qquad
D=\left(\begin{array}{cc}
\lambda_1 & 0\\ 0 & \lambda_2\end{array}\right).$$
%$$\wtF(z)=U(z) P z^D P^{-1} Y_0.$$
By standard matrix theory, we can write $P^{-1}Y_0=TH$ where
$T$ is upper triangular and $H\in SU(2)$. Then \eqref{eq-solution-wtF} gives
\begin{equation}
\label{eq-solution-F}
F(z)=N(z)U(z) P z^D P^{-1} Y_0=N(z) U(z) P z^DTH
\end{equation}
Assume that $\frac{1}{c}\leq |\alpha^2\beta|\leq c$ and $||F(z)||\leq c$ on the unit
circle, for some real number $c$. We need a uniform bound (depending only
on $c$) of $U(z)$ in the unit disk.
The theory of fuchsian systems gives us a bound
of $U(z)$ by construction, but this bound is not uniform as $\alpha\beta\to 0$
(because we are approaching the resonant case).
To obtain a uniform bound, we must use the fact that the monodromy of $F$ is in 
$SU(2)$.
\medskip

First of all, our hypothesis imply the following bounds:
\begin{equation}
\label{eq-bounds}
\frac{1}{2}\leq \Delta\leq \frac{3}{2},\qquad
\frac{-3}{2}\leq \lambda_2\leq \frac{-1}{2},\qquad
|\alpha^{-1}\lambda_1|\leq c\quad\mbox{ and }\quad
\frac{1}{2c}\leq|\alpha\lambda_1|\leq \frac{3c}{2}.
\end{equation}
Hence the matrix $P$ is uniformly bounded.
From \eqref{eq-solution-F}, we obtain
$$1=\det F(z)=z \det (U(z)) z^{-1}\det T.$$
Hence $\det U(z)$ is constant, and since $U(0)=I_2$, we obtain
$\det U(z)=\det T=1.$
The monodromy of $F$ is given by
$$M_{\gamma}(F)=H^{-1} T^{-1}e^{2\pi \i D} TH
=H^{-1}\left(\begin{array}{cc}
e^{2\pi \i \lambda_1} & T_{12}T_{22} (e^{2\pi \i\lambda_1}-e^{2\pi \i\lambda_2})\\
0 & e^{2\pi \i\lambda_2}\end{array}\right)H\in SU(2).$$
Since $\lambda_1-\lambda_2$ is not an integer,
this implies that $T_{12}=0$, so the matrix $T$
is diagonal.
Then $T$ and $z^D$ commute.
Equation \eqref{eq-solution-F} implies that $U(z)PT$ is uniformly bounded
on the unit circle. Since $U(z)PT$ is holomorphic, it is uniformly bounded in the unit disk
by the maximum principle.
Taking $z=0$, we obtain that $T$ is uniformly bounded,
hence $U(z)$ is uniformly
bounded in the unit disk.
Expanding the product in \eqref{eq-solution-F}, we obtain
\begin{equation}
\label{eq-approx-F}
F(z)=\frac{1}{\Delta^{1/4}}\left(\begin{array}{cc}
z^{\lambda_1}(-T_{11}\lambda_2+O(z)) & 
z^{\lambda_1}(T_{22}\alpha\lambda_1 +O(z))\\
z^{1+\lambda_1}(T_{11}\alpha^{-1}\lambda_2+O(z)) &
z^{1+\lambda_2}(-T_{22}\lambda_2 +O(z))\end{array}\right)H.
\end{equation}
where $O(z)$ is holomorphic and uniformly bounded.
Using \eqref{eq-half-space} and the bounds \eqref{eq-bounds}, we obtain
$$(x_1+\i x_2)(z)=-\frac{1}{z}\frac{\alpha\lambda_1}{\lambda_2}(1+O(z)+|z|^{\Delta}O(1))$$
$$x_3(z)=\frac{1}{|z|^{2+2\lambda_2}}\frac{\sqrt{\Delta}}{|T_{22}\lambda_2|^2}
(1+O(z)+|z|^{\Delta}O(1))$$
where $O(z)$ and $O(1)$ are real analytic functions that have uniformly bounded derivatives and $O(z)$ vanishes at the origin.
The conclusions of Theorem \ref{theorem-end} follow.

\medskip

%\subsubsection{Example: catenoid cousins}
%Let us illustrate the above computation in the case of the catenoid cousin family
%(Example 2 in Bryant \cite{bryant}).
%The catenoid cousins are given by
%$$F(z)=\frac{1}{\sqrt{2\mu+1}}\left(\begin{array}{cc}
%(\mu+1)z^{\mu} & \mu z^{-\mu-1}\\
%\mu z^{\mu+1} & (\mu+1) z^{-\mu}\end{array}\right)$$
%where $\mu\in(-\frac{1}{2},\infty)$, $\mu\neq 0$ is a parameter.
%The holomorphic data is
%$$G=\frac{1}{z},\qquad \omega=\mu (\mu+1) dz.$$
%The eigenvalues are
%$$\lambda^-=-1-\mu,\qquad\lambda^+=\mu.$$
%Equation \eqref{eq-solution-F} is satisfied by taking
%$$U(z)\equiv I_2,\qquad T=H=I_2.$$
%The catenoid cousins have two embedded ends, although they are globally embedded only
%for $-\frac{1}{2}<\mu<0$.

\end{document}